\documentclass[10pt]{article}
\usepackage{amsmath,amsfonts,amssymb,amsthm,amscd,enumerate,wasysym,latexsym,amscd,mathrsfs}
\usepackage[dvipsnames]{xcolor}
\usepackage{hyperref}

\usepackage[english]{babel}

\usepackage{graphicx}%
\usepackage{multirow}%
\usepackage{amsmath,amssymb,amsfonts}%
\usepackage{amsthm}%
\usepackage{mathrsfs}%
\usepackage[title]{appendix}%
\usepackage{xcolor}%
\usepackage{textcomp}%
\usepackage{manyfoot}%
\usepackage{booktabs}%
\usepackage{algorithm}%
\usepackage{algorithmicx}%
\usepackage{algpseudocode}%
\usepackage{listings}%

\theoremstyle{definition}
\newtheorem{theorem}{Theorem}%
\newtheorem{proposition}[theorem]{Proposition}%
\newtheorem{example}{Example}%
\newtheorem{remark}{Remark}%
\newtheorem{lemma}{Lemma}
\newtheorem{definition}{Definition}%

\title{Pseudodifferential Operators on Noncommutative Tori: a Survey}

\author{Carolina Neira Jim\'enez\footnote{cneiraj@unal.edu.co, Departamento de Matem\'aticas, Universidad Nacional de Colombia, Carrera 30 \# 45--03, 111321, Bogot\'a, Colombia}}

\date{}

\begin{document}

\maketitle

{\bf Abstract:} {This article presents a survey of recent developments on pseudodifferential operators on noncommutative tori. We describe currently available constructions of those operators: by means of a $C^*$--dynamical system, by using an analogue of the Fourier series representation of a function in the (commutative) torus $C^\infty(\mathbb{T}^n)$, as Rieffel deformations of the standard pseudodifferential operators on $C^\infty(\mathbb{T}^n)$, and in association to certain spectral triples. } \\



\section{Introduction}\label{sec1}

Some of the most prototypical examples of noncommutative spaces are given by noncommutative tori, which have been widely studied (see e.g.\ \cite{cohenconnes}, \cite{connesmoscovicimodularcurvature}, \cite{essouabriiochumlevysitarz}, \cite{fathizadehkhalkhaligb}, \cite{fathizadehkhalkhalisc2t}, \cite{haleepongeNCTI}, \cite{haleepongeNCTII}, \cite{levyneirapaycha}, \cite{liu}, \cite{plazas}, \cite{rieffelnctori}, \cite{ponge} and references therein). They can be defined as deformations of the algebra of smooth functions on the usual torus, and they can be seen as the underlying noncommutative algebra of certain spectral triples. \\

As discussed in  \cite[p.\ 523]{graciabondia}, the term noncommutative torus refers to two species of algebras: the $C^*$--algebras originally called irrational rotation algebras, as well as the pre--$C^*$--algebras that consist of the smooth elements of these $C^*$--algebras under the action of certain Lie groups of automorphisms. The second one is used when we focus on the differential structure. \\

Among the constructions of noncommutative tori, we find a deformation in the sense of Rieffel. In his monograph \cite{rieffel93}, M.\ Rieffel describes the construction of strict deformation quantizations in the setting in which the Poisson bracket on a manifold is defined by a smooth action of the Euclidean space on the manifold, and he shows that this construction can be also carried out for actions of the Euclidean space on any $C^*$--algebra. For the construction of noncommutative tori, this deformation mainly consists in inserting an appropriate twist of the (commutative) pointwise product in the algebra of smooth functions on the torus $C^\infty(\mathbb{T}^n)$, which extends it to a new associative and noncommutative product. \\

Additionally, it is common to look for noncommutative counterparts of known concepts, such as pseudodifferential operators. On an open subset $U$ of $\mathbb{R}^n$, a symbol of order $d\in\mathbb{R}$ is a smooth function $\sigma : T^*U \to\mathbb{C}$, such that for all $\alpha, \beta\in\mathbb{N}^n$, and for every compact subset $K\subset U$, there exists $C_{\alpha,\beta,K}\in \mathbb{R}$ such that 
\begin{equation}\label{eq:symbolonopenset}
 \left\lvert\partial_x^\alpha\partial_\xi^\beta \sigma (x,\xi)\right\rvert \leq C_{\alpha,\beta,K} (1+\lvert\xi\rvert)^{d-\lvert\beta\rvert},\quad\text{ for all } x\in K,\ \xi\in T_x^*U\cong\mathbb{R}^n.
\end{equation}
The pseudodifferential operator $P_\sigma$ associated to the symbol $\sigma$ is given by
\begin{equation}\label{eq:def pdos Rn}
P_\sigma(f)(x):=\int_{\mathbb{R}^n}\int_U e^{-2\pi i\langle x-y,\xi\rangle}\sigma(x,\xi)f(y)\,dy\,d\xi,
\end{equation}
where $f\in C_0^\infty(U)$, $x\in U$, and we consider a normalized volume measure on the Euclidean space. On a closed manifold, a pseudodifferential operator is a linear operator that locally can be written as \eqref{eq:def pdos Rn} (see e.g.\ \cite{shubin}). In the literature these operators are known as {standard pseudodifferential operators} on $\mathbb{R}^n$. \\

The description of the manifold by local charts is crucial in this definition of the symbols of pseudodifferential operators. As there is no underlying geometric space on a noncommutative space like a noncommutative torus, such a representation by local charts is not available. There are nevertheless alternative approaches to pseudodifferential calculi on that space: by means of a $C^*$--dynamical system (\cite{connestretkoff}), by using an analogue of the Fourier series representation of a function in the (commutative) torus $C^\infty(\mathbb{T}^n)$ (\cite{levyneirapaycha}), as a Rieffel deformation of the standard pseudodifferential operators on $C^\infty(\mathbb{T}^n)$ (\cite{liu}), and in association to a spectral triple (\cite{connesmoscovici}). In this document we revisit those definitions, in order to provide an overview of various pseudodifferential techniques that can be defined on certain noncommutative spaces. \\

The document is structured as follows: In the second section we revisit Rieffel's deformation in order to see noncommutative tori as a particular case. In the third section we give the definition of noncommutative tori as deformations of the algebra of continuous functions on the commutative torus starting both from $\mathbb{T}^n$ and from its Pontryagin dual $\mathbb{Z}^n$. In Section 4 we revisit various calculi of pseudodifferential operators on noncommutative tori. In Section 5 we consider noncommutative tori as the underlying noncommutative algebra of certain spectral triples, and finally in Section 6, we recall the definition of abstract pseudodifferential operators on regular spectral triples.

\section{Rieffel deformation}\label{sect:rieffel deformation}

As mentioned in the Introduction, the deformation described by M.\ Rieffel in \cite{rieffel93} mainly consists in inserting an appropriate twist of the product in a given algebra, in order to extend it to a new associative and noncommutative product. It is important to emphasize that this twist of the product is defined from the action of a vector space, and later in Section \ref{sect:NCT as a deformation} we will consider a deformation defined from the action of the torus as an abelian group. In the present section we revisit the deformation defined from the action of a vector space, both in a differentiable Fr\'echet algebra and in a $C^*$--algebra, in order to see both noncommutative tori (Section \ref{sect:NCT as a deformation} below) and a pseudodifferential calculus on those spaces (Section \ref{sect:deformed pdos} below) as particular cases. We closely follow \cite{rieffel93}. 

\subsection{Deformation of a differentiable Fr\'echet algebra along the action of a vector space}\label{sect:Afrechet}

Let $V$ be a real vector space with basis $\{X_1,\ldots, X_n\}$; we will also consider $V$ as a Lie group. Let $A$ be a Fr\'echet algebra over $\mathbb{C}$ endowed with a family of seminorms $\{\vert\ \vert_j\}_{j\in\mathbb{N}}$ which determine the topology of $A$. In particular, the product is jointly continuous, so that for every $j\in\mathbb{N}$ there exist $k\in\mathbb{N}$ and a constant $c_j\in\mathbb{R}$ such that for all $a,b\in A$
\begin{equation*}
\vert ab\vert_j\leq c_j\vert a\vert_k\vert b\vert_k.
\end{equation*}
Let 
\begin{align*}
\alpha:V\times A&\to A \\
(v,a)&\mapsto\alpha_v(a)
\end{align*}
be an isometric action of $V$ on $A$ as automorphisms of $A$, and denote by $A^\infty$ the space of smooth (meaning $C^\infty$) elements for the action $\alpha$ of $V$ on $A$:
\begin{equation}\label{def:Ainfty}
A^\infty:=\{a\in A:v\mapsto \alpha_v(a) \text{ is smooth}\}.
\end{equation}

For each $1\leq k\leq n$, let $\partial_k:=\alpha_{X_k}$ denote the usual operator of \emph{partial differentiation} on $A^\infty$ in the direction of $X_k$, (where we identify the Lie group $V$ with its Lie algebra). For $\mu = (\mu_1,\ldots,\mu_n)\in\mathbb{N}^n$ we let $\partial^\mu:=\partial_1^{\mu_1}\partial_2^{\mu_2}\cdots\partial_n^{\mu_n}$ denote the usual corresponding higher partial differentiation. 
Then, the action $\alpha$ is strongly continuous, namely, for all $a\in A^\infty$, given $\mu\in\mathbb{Z}^n$, $j\in\mathbb{N}$, there exist $m\in\mathbb{N}$ and a constant $C_{\mu,j,m}$ such that
\begin{equation*}
\left\vert\partial^\mu\alpha_v(a)\right\vert_j\leq C_{\mu,j,m}\vert a\vert_{m}, \text{ for all }v\in V.
\end{equation*}

We equip $A^\infty$ with the seminorms \cite[Chapt.\ 1]{rieffel93}
\begin{equation*}
\vert a\vert_{j,k}:=\sup_{i\leq j}\sum_{\lvert\mu\rvert\leq k}\dfrac{1}{\mu!}\vert\partial^\mu a\vert_i,
\end{equation*}
which endow $A^\infty$ with a Fr\'echet algebra structure. Moreover, the action $\alpha$ is invariant on $A^\infty$, it turns out to be isometric for these seminorms, and it is not only strongly continuous, but also differentiable in the sense that each element of $A^\infty$ is a smooth element for this action and these seminorms. \\ 

We denote by $C_b(V,A)$ \emph{the Fr\'echet space of continuous bounded functions from $V$ to $A$}, equipped with the seminorms $\vert f\vert_k:=\sup\limits_{v\in V}\vert f(v)\vert_k$, for $k\in\mathbb{N}$, where, on the right hand side we have the corresponding seminorm on $A$. Let $\tau$ denote the action of $V$ on $C_b(V, A)$ given by translation:
\begin{equation*}
(\tau_tf)(x):=f(x+t),\qquad t,x\in V.
\end{equation*}
This action is isometric for each of the above seminorms, but it is not strongly continuous.  
We then consider $C_u(V,A)$, the largest subspace of $C_b(V,A)$ on which $\tau$ is strongly continuous, i.e.\ the space of \emph{$A$--valued uniformly continuous bounded functions on $V$}, and consider $\tau$ as an action of $V$ on $C_u(V,A)$. The space $\mathcal{B}^A(V)$ denotes the Fr\'echet subalgebra of $C_u(V,A)$ consisting of smooth elements for the action $\tau$, and applying the above considerations to $W:=V\times V$ we consider $\mathcal{B}^A(W)$. \\

We fix an inner product $\langle\ ,\ \rangle$ on $V$. By choosing an appropriate Haar measure on $V$, and hence on $W$, for any $F\in\mathcal{B}^A(W)$ the integral 
\begin{equation*}
\int_W F(u,v)\,e^{2\pi i\langle u, v\rangle}\,du\,dv
\end{equation*}
is well defined as an oscillatory integral \cite[Def.\ 1.3]{rieffel93}. If $J$ is a linear operator on $V$, for any $a,b\in A$, the function $F(u,v):=\alpha_{Ju}(a)\alpha_v(b)$ belongs to $\mathcal{B}^A(W)$, and thus, the integral 
\begin{equation*}
\int_V\int_V \alpha_{Ju}(a)\alpha_v(b)\,e^{2\pi i\langle u, v\rangle}\,du\,dv\in A 
\end{equation*}
is well defined \cite[Chapt.\ 2]{rieffel93}. \\

We are now ready to define a deformed product between two elements of the Fr\'echet algebra $A$ with action $\alpha$ from the Lie group $V$. \\

\begin{definition}[{\cite[Def.\ 2.1]{rieffel93}}]\label{def:deformedproduct}
Let $a,b\in A$ and let $J$ be a linear operator on $V$. The \emph{deformed product} of $a$ and $b$ is defined by
\begin{align*}
a\times_J b&:=\int_W \alpha_{Ju}(a)\alpha_v(b)\,e^{2\pi i\langle u, v\rangle}\,du\,dv\\
&=\int_W \alpha_{u}(a)\alpha_{J^tv}(b)\,e^{2\pi i\langle u, v\rangle}\,du\,dv.
\end{align*}
\end{definition} 

The product $\times_J$ is bilinear, continuous for the given Fr\'echet topology on $A$, and associative. $V$ acts (differentially) on $A$ by automorphisms of the deformed product, and if $J=0$, we get back the product of the algebra $A$: $a\times_0b=ab$. \\

If $J$ is skew symmetric and if $A$ has a continuous involution, the differentiable isometric action $\alpha$ of $V$ on $A$ is by $*$--automorphisms. In this case, the involution on $A$ is also an involution for the deformed product $\times_J$. In the following, we will consider that $J$ is skew symmetric.\\

If $L$ is an integer lattice for $V=\mathbb{R}^n$, then $H:=V/L$ is a compact group, with dual group $L$. If $L\subseteq\ker(\alpha)$, $\alpha$ can be viewed as an action of $H$. For $k\in L$, the corresponding character of $H$ is the function $v\mapsto e^{2\pi i\langle k, v\rangle}$, and the associated spectral subspace is given by the set of eigenvectors of this action corresponding to the eigenvalue $e^{2\pi i\langle k, v\rangle}$, also known as the \emph{$k$--th isotypic component} of $A$ \cite[Chapt.\ 2]{rieffel93}:
\begin{equation}\label{eq:spectralsubspace}
A_k:=\{a\in A:\alpha_v(a)=e^{2\pi i\langle k, v\rangle}a \text{ for all }v\in V\}.
\end{equation}
The direct sum of all the $A_k$'s is dense in the algebra $A$. Moreover, for $k,l\in L$, if $a\in A_k$ and $b\in A_l$ then 
\begin{equation}\label{eq:product of eigenvectors}
a\times_Jb=e^{-2\pi i\langle k, Jl\rangle}ab.
\end{equation}

For $f,g\in\mathcal{B}^A(V)$ we obtain a well defined product $f\times_J g$ in $\mathcal{B}^A(V)$:
\begin{align*}
(f\times_Jg)(x)&:=\iint \tau_{Ju}(f)(x)\tau_v(g)(x)e^{2\pi i\langle u, v\rangle}\,du\,dv\\
&=\iint f(x+Ju)g(x+v)e^{2\pi i\langle u, v\rangle}\,du\,dv.
\end{align*}
With this deformed product and with its usual seminorms, $\mathcal{B}^A(V)$ becomes a Fr\'echet algebra. Moreover, the action $\tau$ is by isometric algebra automorphisms for this deformed product. We denote by $\mathcal{S}^A(V)$ the subspace of $\mathcal{B}^A(V)$ consisting of functions such that the products of their derivatives with any polynomial on $V$ are bounded for each of the seminorms on $A$ \cite[Chapt.\ 3]{rieffel93}. \\

\subsection{Deformation of a $C^*$--algebra along the action of a vector space}
\label{sect:deformation cstar algebra}

Let $A$ be a $C^*$--algebra, let $\alpha$ be an action of $V$ on $A$ as a Lie group, and let $J$ be a skew symmetric linear operator on $V$. The goal now is to realize $\mathcal{B}^A(V)$ as an algebra of bounded operators on a Hilbert space, and endow it with a $C^*$--norm. Consider the left action of $\mathcal{B}^A(V)$ on $\mathcal{S}^A(V)$ given by the deformed product: $L_F^J(g):=F\times_J g$, where $F\in\mathcal{B}^A(V)$, $g\in\mathcal{S}^A(V)$, and consider the map \cite[Chapt.\ 4]{rieffel93}
\begin{align}
A&\to C_u(V,A) \notag\\
a&\mapsto\tilde{a}, \label{def:maptilde}
\end{align}
where to any $a\in A$ we associate the function $\tilde{a}(x):=\alpha_x(a)$ on $V$. \\

The map \eqref{def:maptilde} is a $*$--homomorphism which is equivariant for the action $\alpha$ on $A$ and the action $\tau$ on $C_u(V,A)$. Thus, it carries smooth vectors to smooth vectors ($A^\infty$ into $\mathcal{B}^A(V)$), and is a homomorphism for their deformed products $\times_J$. By this correspondence, each $a\in A^\infty$ determines a bounded operator $L_a^J$ on $\mathcal{S}^A(V)$, by
\begin{align*}
L_a^J(f)(x)&:=(\tilde{a}\times_Jf)(x) \qquad  \qquad  \qquad \text{ for } f\in\mathcal{S}^A(V), \ x\in V,\\
&=\iint \tilde{a}(x+Ju)f(x+v)e^{2\pi i\langle u, v\rangle}\,du\,dv\\
&=\iint \alpha_{x+Ju}(a)\tau_v(f)(x)e^{2\pi i\langle u, v\rangle}\,du\,dv.
\end{align*}
We can define a new norm $\|\ \|_J$ on $A^\infty$ by $\|a\|_J:=\|L_a^J\|$ (operator norm), so that $A^\infty$ with the product $\times_J$ and this new norm is a pre--$C^*$--algebra. \\

\begin{definition}[{\cite[Def.\ 4.9]{rieffel93}}] 
The completion of $A^\infty$ for the norm $\|\ \|_J$, denoted by $A_J$, is called the \emph{deformation} of $A$ by $J$ (and $\alpha$), or the \emph{algebra obtained from $A$ by deformation quantization}. \\
\end{definition}

The algebra $A_J$ can be identified via the map \eqref{def:maptilde} with a $C^*$--subalgebra of the completion of $\mathcal{B}^A_J(V)$ for the operator norm. \\

\begin{remark}
In this section we have followed the approach given in \cite{rieffel93} which mainly deals with the analytical difficulties when the action that produces the deformation is along the Euclidean space. When the action is along a compact group, like the torus, we can work with its Pontryagin dual given by the additive group $\mathbb{Z}^n$ as we do in the following sections. In that case, we work with a Fourier series representation of an element of the algebra.
\end{remark}

\section{Noncommutative tori}\label{sect:NCT as a deformation}

Some of the most prototypical examples of noncommutative spaces are given by noncommutative tori, which have been widely studied. We refer the reader for example to \cite{connesncg}, \cite{graciabondia} and \cite{rieffelnctori}. In \cite{haleepongeNCTI} and \cite{ponge} there is a good summary of references about different studies of noncommutative tori both in mathematics and in mathematical physics. In this section we consider those spaces as deformations of the algebra of continuous functions on the commutative torus, following the construction of the previous section. However, as mentioned before, here we will consider a deformation defined from the action of the torus as an abelian compact group. \\

In this paper we present two constructions of deformation: $\theta$--deformation and 2--cocycle deformation for noncommutative tori. The $\theta$--deformation is a deformation along an action of an abelian group, usually a torus action. With the action of a torus, we can apply the $\theta$--deformation on the coordinate functions using a 2--cocycle to deform the multiplication of the algebra as we explained in the previous section. We refer the reader to \cite{wilson} for a review of compact quantum groups that arise as $\theta$--deformations of compact Lie groups such as the Connes--Landi spheres \cite{conneslandi}.  \\

In noncommutative geometry, thanks to the Gelfand--Naimark theorem, noncommutative $C^*$--algebras can be regarded as formal duals to spaces. If $A$ is a $C^*$--algebra, and $\widehat{A}$ denotes the set of unital homomorphisms from $A$ to $\mathbb{C}$, then $\widehat{A}$ is a compact Hausdorff space. The Gelfand--Naimark theorem states that for every commutative $C^*$--algebra $A$, there is an isometric $*$--isomorphism between $A$ and the algebra of complex valued continuous functions on $\widehat{A}$, denoted by $C(\widehat{A})$ (see e.g.\ \cite[Thm.\ D.5.11]{ruzhanskyturunen}). \\

Thus, by Gelfand--Naimark duality, as a topological space, the $n$--dimensional torus $\mathbb{T}^n:=\mathbb{R}^n/\mathbb{Z}^n$ is a compact Hausdorff space characterized by $C(\mathbb{T}^n)$, the unital commutative $C^*$--algebra of continuous complex valued functions on $\mathbb{T}^n$. \\  

Since the Pontryagin dual of $\mathbb{T}^n$ is isomorphic to $\mathbb{Z}^n$ (see e.g.\ \cite[Rem.\ 7.5.9]{ruzhanskyturunen}), we will present the definition of noncommutative tori as deformations of the algebra of continuous functions on the commutative torus $C(\mathbb{T}^n)$ starting both from $\mathbb{T}^n$ and from $\mathbb{Z}^n$.

\subsection{Deformation of functions on $\mathbb{T}^n$ along the action of $\mathbb{R}^n$}
\label{section:NCTalg}

With the notation of the previous section, consider $V=\mathbb{R}^n$ and the action $\alpha$ given by translation on the torus $\mathbb{T}^n$. If $A=C(\mathbb{T}^n)$, as in \eqref{def:Ainfty} the subalgebra of smooth elements of $A$ for the action of $\mathbb{T}^n$ is $A^\infty=C^\infty(\mathbb{T}^n)$, and the action $\alpha$ factors through the compact group $H=\mathbb{T}^n=\mathbb{R}^n/\mathbb{Z}^n$.  \\

For $k\in \mathbb{Z}^n$, let $E_k:\mathbb{T}^n\to\mathbb{C}$ be given by $E_k(x):=e^{2\pi i\langle x,k\rangle}$. As in \cite[Ex.\ 3.5]{liu}, by elementary Fourier theory on $\mathbb{T}^n$, the set $\{E_k\}_{k\in\mathbb{Z}^n}$ serves as a basis for $C(\mathbb{T}^n)$, that is, if $f\in C(\mathbb{T}^n)$, there exists a sequence of complex numbers $(f_k)_{k\in\mathbb{Z}^n}$ such that the Fourier series
\begin{equation}\label{eq:decomposition f torus}
\sum\limits_{k\in\mathbb{Z}^n}f_{k}E_k
\end{equation}
converges to the function $f$. Moreover, if $f\in C^\infty(\mathbb{T}^n)$, the sequence $(f_k)_{k\in\mathbb{Z}^n}$ belongs to $\mathcal{S}(\mathbb{Z}^n)$, the vector space of sequences $(f_k)_{k\in\mathbb{Z}^n}$ that decay faster than the inverse of any polynomial in $k$. \\

\begin{remark}[{\cite[Rem.\ 7.5.9]{ruzhanskyturunen}}]
The functions $E_k$, for $k\in\mathbb{Z}^n$, are used to identify the Pontryagin dual of $\mathbb{T}^n$ with $\mathbb{Z}^n$. \\
\end{remark}

The action $\alpha$ on the elements $E_k$ is given by
\begin{align}
(\alpha_x(E_k))(y)&=E_k(y+x)=e^{2\pi i\langle y+x, k\rangle} \notag\\
&=e^{2\pi i\langle x, k\rangle}e^{2\pi i\langle y, k\rangle}=e^{2\pi i\langle x, k\rangle}E_k(y), \label{eq:action torus Ek}
\end{align}
for $x,y\in\mathbb{T}^n$, $k\in\mathbb{Z}^n$. From this we recognize that $\alpha$ coincides with the usual dual action of the dual group $\mathbb{T}^n$ of $\mathbb{Z}^n$ on $C(\mathbb{T}^n)$, except that it is being viewed as an action of $\mathbb{R}^n$ instead of $\mathbb{T}^n$.\\

Following the notation of Section \ref{sect:deformation cstar algebra}, let $J:={\theta}/{2}$, where $\theta$ is an $n\times n$ skew symmetric real matrix. For $k,l\in \mathbb{Z}^n$, consider the deformed product
\begin{equation}\label{eq:EkEl}
E_k\times_\theta E_l:=e^{-2\pi i\langle k, Jl\rangle}E_{k+l}=e^{-\pi i\langle k,\theta l\rangle}E_{k+l}.
\end{equation}
For all $k\in \mathbb{Z}^n$, $E_k$ is a unitary element of the deformed algebra $A_\theta$, which is isomorphic to the deformed algebra $A_J$ \cite[Ex.\ 10.2]{rieffel93}. \\

\begin{definition}[{\cite[Ex.\ 10.2]{rieffel93}}]\label{def:quantum torus}
Given $n\in\mathbb{N}$ and an $n\times n$ skew symmetric real matrix $\theta$, the \emph{noncommutative $n$--torus with parameter $\theta$}, or more precisely, the algebra of continuous functions on the noncommutative $n$--torus with parameter $\theta$, is the universal $C^*$--algebra $A_\theta^n$ generated by unitaries satisfying the relation given by \eqref{eq:EkEl}. The subalgebra of smooth elements of $A_\theta^n$ for the action of $\mathbb{T}^n$ given in \eqref{def:Ainfty} and denoted by $\mathcal{A}_\theta^n$, is known as the \emph{smooth noncommutative $n$--torus with parameter $\theta$}, or more precisely, the algebra of smooth functions on the noncommutative $n$--torus with parameter $\theta$. \\
\end{definition} 

The elements $\{E_k\}_{k\in\mathbb{Z}^n}$ are examples of the generators mentioned in Definition \ref{def:quantum torus}.  \\

\begin{remark}\label{rem:algebrafunctns}
In this definition it is important to notice that, as there is no underlying geometric space, the use of the terms ``algebra of continuous functions”  and ``algebra of smooth functions” is only metaphorically speaking, and we actually refer to the space of elements \emph{playing the role} of continuous resp.\ smooth functions on the noncommutative space. \\
\end{remark}

\begin{remark}
Another important point is that Definition \ref{def:quantum torus} gives actually a family of spaces parametrized by the dimension $n$ and the deformation parameter $\theta$, and usually it is more precise to write \emph{noncommutative tori}. In the following we take this point into account fixing the dimension and the parameter.  \\
\end{remark}

Notice that, by Gelfand--Naimark duality, $A_0^n=C(\mathbb{T}^n)$ corresponds to the usual algebra of continuous functions on the commutative torus. If $n=2$ and the entries of $\theta$ are irrational, the family $A_\theta^n$ of noncommutative tori are also known as irrational rotation algebras (see \cite{cohenconnes}, \cite{pimsnervoiculescu}, or \cite{rieffelirrationalrot}). \\ 

The deformed algebra $\mathcal{A}_\theta^n$ is identical to $C^\infty(\mathbb{T}^n)$ as a topological vector space with the deformed product
\begin{align}
f\times_\theta g &:=\sum_{k,l\in\mathbb{Z}^n}f_kg_lE_k\times_\theta E_l \notag \\
&=\sum_{k,l\in\mathbb{Z}^n}e^{-\pi i\langle k,\theta l\rangle}f_kg_lE_{k+l}, \label{eq:product times theta}
\end{align}
where $f_k$ and $g_l$ are the Fourier coefficients of $f$ and $g$, respectively. If we take $f=E_k$ and $g=E_l$, where $k,l\in \mathbb{Z}^n$, we obtain that
\begin{equation}\label{eq:EkEl commuting relation}
E_k\times_\theta  E_l = e^{-2 \pi i\,\langle k,\theta l\rangle}E_l\times_\theta E_k.
\end{equation}

\subsection{A 2--cocycle deformation of functions on $\mathbb{Z}^n$}
\label{section:NCTdiff}

Consider the space $C_c(\mathbb{Z}^n):=\{a:\mathbb{Z}^n\to\mathbb{C}\text{ with finite support}\}$ equipped with the $L^1$--norm
\begin{equation*}
\|a\|_1:=\sum_{k\in\mathbb{Z}^n}\lvert a(k)\rvert.
\end{equation*}
Let $L^1(\mathbb{Z}^n)$ be the involutive Banach algebra defined by the completion of $C_c(\mathbb{Z}^n)$ with respect to this norm.\\

The supremum over all involutive representations $\rho$ of $L^1(\mathbb{Z}^n)$ defines a seminorm on $L^1(\mathbb{Z}^n)$:
\begin{equation*}
\|b\|:=\sup_\rho\lvert\rho(b)\rvert. 
\end{equation*}
In many cases this is already a norm, and the completion of $L^1(\mathbb{Z}^n)$ with respect to that norm is a $C^*$--algebra called \emph{the enveloping $C^*$--algebra of $L^1(\mathbb{Z}^n)$} \cite[Def.\ 12.2]{graciabondia}. \\

Let $\theta$ be an $n\times n$ skew symmetric real matrix. Consider the 2--cocycle
\begin{equation}\label{def:cocycle}
  c(k,l) := e^{- \pi i\,\langle k,\theta l\rangle},\quad\text{ for all } k,l\in \mathbb{Z}^n .
\end{equation}
This is a function $c:\mathbb{Z}^n\times\mathbb{Z}^n\to S^1$ such that $c(0,0)=1$ and $c(r,s)c(r+s,t)=c(r,s+t)c(s,t)$ for all $r,s,t\in\mathbb{Z}^n$. This implies that for all $t\in\mathbb{Z}^n$, $c(0,t)=c(t,0)=1$ and $c(t,-t)=c(-t,t)$. Moreover, the cocycle given in \eqref{def:cocycle} satisfies that for all $r,s\in\mathbb{Z}^n$, $\overline{c(r,s)}=c(s,r)$ and $c(r,-r)=1$. \\

The Banach algebra $L^1(\mathbb{Z}^n,c)$ is defined by introducing a twisted convolution and involution \cite[Def.\ 12.4]{graciabondia}
\begin{align}
(a\star b)(k)&:=\sum_{l\in\mathbb{Z}^n}a(l)b(k-l)c(l,k-l) \notag\\
&=\sum_{l\in\mathbb{Z}^n}a(l)b(k-l)e^{- \pi i\,\langle l,\theta (k-l)\rangle}, \label{twistedconvolution}\\
a^*(k)&:=\overline{c(k,-k)a(-k)} \notag \\
&=\overline{a(-k)}, \label{twistedinvolution}
\end{align}
for any $k\in\mathbb{Z}^n$. Let $C^*(\mathbb{Z}^n,c)$ denote the enveloping $C^*$--algebra of $L^1(\mathbb{Z}^n,c)$. For each $k\in\mathbb{Z}^n$ the function $x\mapsto e^{2\pi i\langle x, k\rangle}$ will correspond to a unitary element $U_k$ in $C^*(\mathbb{Z}^n,c)$. \\

If $\{U_{e_1},\ldots,U_{e_n}\}$ are unitary elements corresponding to the standard basis $\{e_i\}_{1\leq i\leq n}$ for $\mathbb{Z}^n$, they generate $C^*(\mathbb{Z}^n,c)$: If $k\in\mathbb{Z}^n$, denote by $e(k,\theta)$ the function $\exp\left({-\pi i\sum\limits_{\substack{j,l=1 \\ j<l}}^nk_j\theta_{jl}k_l}\right)$; then \cite[p.\ 534]{graciabondia}
\begin{equation}\label{eq:weyl elements}
U_k=U_{(k_1,\ldots,k_n)}:=e(k,\theta)U_{e_1}^{k_1}\cdots U_{e_n}^{k_n},
\end{equation}
define the \emph{Weyl elements} which satisfy that $U_0=1$, and for all $k,l\in\mathbb{Z}^n$
\begin{equation}
\ U_k U_{l} = c(k,l)U_{k+l},
\end{equation}
(compare with Equation \eqref{eq:EkEl}), which, as in \eqref{eq:EkEl commuting relation}, implies that 
\begin{equation}\label{eq:product weyl elements}
U_k U_{l} = e^{-2 \pi i\,\langle k,\theta l\rangle}U_lU_k.
\end{equation}

\begin{remark}
In references such as \cite{haleepongeNCTI}, \cite{haleepongeNCTII}, \cite{tao}, the function $e(k,\theta)$ does not appear in Equation \eqref{eq:weyl elements}. \\
\end{remark}

\begin{remark}
If we set $\omega(k,l):=2\pi\langle k,\theta l\rangle$, according to \cite[Def.\ 4.2.1]{baerginouxpfaffle} we have that a Weyl system for the symplectic vector space $(\mathbb{Z}^n,\omega)$ consists of the  $C^*$--algebra $C^*(\mathbb{Z}^n,c)$ and the map $k\mapsto U_k$ for $k\in\mathbb{Z}^n$. Moreover, $C^*(\mathbb{Z}^n,c)$ becomes a CCR--representation for $(\mathbb{Z}^n,\omega)$ (see \cite[Def.\ 4.2.8]{baerginouxpfaffle} and \cite[Sect.\ 5.2.2.2]{brattellirobinson}). \\
\end{remark}

For $f\in C^\infty(\mathbb{T}^n)$, consider the Fourier transform
\begin{equation}\label{eq:Fourier transform Tn}
\hat{f}(k)=\int_{\mathbb{T}^n}e^{-2\pi i\langle x, k\rangle}f(x)\,dx,\ \  \text{  for } k\in\mathbb{Z}^n,
\end{equation}
which carries the pointwise multiplication on $C^\infty(\mathbb{T}^n)$ to the convolution on $\mathcal{S}(\mathbb{Z}^n)$.\\ 

For every $h\in\mathbb{R}$, consider the cocycle \cite[Ex.\ 2]{rieffeldefquant}
\begin{equation*}
c_h(k,l) := e^{- \pi ih\, \langle k,\theta l\rangle},\quad\text{ for all } k,l\in \mathbb{Z}^n .
\end{equation*} 
As in Equation \eqref{twistedconvolution}, for $\phi,\psi\in\mathcal{S}(\mathbb{Z}^n)$ we set
\begin{equation}\label{eq:hconvolution product on Zn}
(\phi\star_{h}\psi)(k):=\sum_{l\in\mathbb{Z}^n}\phi(l)\psi(k-l)c_h(l,k-l),\qquad k\in\mathbb{Z}^n.
\end{equation}
We define the involution on $\mathcal{S}(\mathbb{Z}^n)$ to be that coming from complex conjugation on $C^\infty(\mathbb{T}^n)$, so that, as in Equation \eqref{twistedinvolution}, for $\phi\in\mathcal{S}(\mathbb{Z}^n)$
\begin{equation*}
\phi^*(k)=\overline{\phi}(-k),\qquad k\in\mathbb{Z}^n.
\end{equation*}
Define the norm $\vert\ \vert_h$ on $\mathcal{S}(\mathbb{Z}^n)$ to be the operator norm for the action of $\mathcal{S}(\mathbb{Z}^n)$ on $\ell^2(\mathbb{Z}^n)$ given by Formula \eqref{eq:hconvolution product on Zn}.   \\ 

Let us denote by $C_h$ the space $C^\infty(\mathbb{T}^n)$ equipped with the product, involution, and norm obtained by pulling back through the inverse Fourier transform the product $\star_h$, involution and norm $\vert\ \vert_h$ described above. The completion of the $C_h$'s form a continuous field of $C^*$--algebras, $C_0$ is just $C^\infty(\mathbb{T}^n)$ with the usual pointwise multiplication and supremum norm, $C_1=\mathcal{A}_\theta^n$ with the notation of Definition \ref{def:quantum torus}, and $A_\theta^n$ becomes its norm completion as a $C^*$--algebra. \\

The $n\times n$ skew symmetric real matrix $\theta$ allows us to define a Poisson bracket between two functions $f,g\in C^\infty(\mathbb{T}^n)$ by \cite[Ex.\ 2]{rieffeldefquant}
\begin{equation}\label{poisson bracket}
\{f,g\}_\theta:=\sum_{j,k=1}^n\theta_{jk}\dfrac{\partial f}{\partial x_j}\dfrac{\partial g}{\partial x_k},
\end{equation}
where $(x_1,\ldots,x_n)$ denotes the local coordinates on $\mathbb{T}^n$. The deformation of the pointwise product of $C^\infty(\mathbb{T}^n)$ to the one--parameter family of associative products $\star_{h}$ satisfies that \cite[Chapt.\ 9]{rieffel93}
\begin{equation*}
\left\vert\dfrac{f\star_hg-g\star_hf}{ih}-\{f,g\}_\theta\right\vert_h\to0, \qquad \text{ as }h\to0.
\end{equation*}
Therefore, the $C_h$'s form a \emph{strict deformation quantization} of $C^\infty(\mathbb{T}^n)$ in the direction of the Poisson bracket \eqref{poisson bracket} \cite[Def.\ 9.2]{rieffel93} (see also \cite{rieffelnctori}). \\

The action of $\mathbb{T}^n=\mathbb{R}^n/\mathbb{Z}^n$ on $C^*(\mathbb{Z}^n,c)$ is given by
\begin{equation}\label{eq:actiontorus}
\alpha_v(U_k):=e^{\,2\pi i \langle v, k\rangle} U_k, \quad\text{ for all } v\in \mathbb{T}^n.
\end{equation}
and as in Definition \ref{def:quantum torus}, the algebra of smooth elements of $C^*(\mathbb{Z}^n,c)$ for this action is the smooth noncommutative torus $\mathcal{A}_\theta^n$, which satisfies
\begin{equation*}
\mathcal{A}_\theta^n=\left\{\sum_{k\in \mathbb{Z}^n} a_k U_k\in L^1(\mathbb{Z}^n,c):\  (a_k)_k\in \mathcal{S}(\mathbb{Z}^n)\right\}.
\end{equation*}

Let us consider the map
\begin{align}
\mathbf{t}:\ L^1(\mathbb{Z}^n,c)\ &\to\mathbb{C} \notag \\
a=\sum_{k\in \mathbb{Z}^n} a_k U_k&\mapsto\mathbf{t}(a) := a_0, \label{eq:trace t}
\end{align}
which extends as a (normalized) trace on $C^*(\mathbb{Z}^n,c)$. With this trace we can define an inner product on $C^*(\mathbb{Z}^n,c)$:
\begin{equation}\label{eq:inner product t}
 \langle a,b\rangle_\mathbf{t}:=\mathbf{t}(ab^*) \quad\text{ for all }  a,b\in C^*(\mathbb{Z}^n,c),
\end{equation}
with associated norm $\|\ \|_\mathbf{t}$. \\

Similarly to the derivatives of the Fourier transform of a function, the infinitesimal generators of the action \eqref{eq:actiontorus} are given by $\{2\pi i\delta_j\}_{1\leq j\leq n}$, where \cite[Sect.\ 2]{levyneirapaycha}
\begin{equation}\label{eq:deltas}
 \delta_j \left(\sum_{k\in \mathbb{Z}^n} a_k U_k\right):= \sum_{k\in \mathbb{Z}^n} a_k k_j  U_k.
\end{equation}
Each $\delta_j$ is a $*$--derivation of $\mathcal{A}_\theta^n$, that is, it satisfies
\begin{enumerate}
\item $(\delta_j(a))^*=\delta_j(a^*)$,
\item  $\delta_j(ab)=\delta_j(a)b+a\delta_j(b)$,
\end{enumerate}
for all $a,b\in\mathcal{A}_\theta^n$. These directional derivatives can be used to define partial differential operators on $\mathcal{A}_\theta^n$.  \\

\begin{example}
\begin{enumerate}
\item A first order differential operator on $\mathcal{A}_\theta^n$ is of the form $\sum\limits_{j=1}^na_j\delta_j$, for $a_j\in \mathcal{A}_\theta^n$.
\item The Laplace operator on $\mathcal{A}_\theta^n$ is $\Delta=\sum\limits_{j=1}^n\delta_j^2$.
\item If $r\in\mathbb{N}$, a differential operator on $\mathcal{A}_\theta^n$ of order at most $r$ is given by $\sum\limits_{\lvert\mu\rvert\leq r}a_\mu\delta_1^{\mu_1}\cdots\delta_n^{\mu_n}$, for $a_\mu\in \mathcal{A}_\theta^n$. \\
\end{enumerate}
\end{example}

Moreover, these derivations satisfy that $\mathbf{t} \circ \delta_j\ = 0$, and we can equip $\mathcal{A}_\theta^n$ with a structure of Fr\'echet $*$--algebra given by the following family of seminorms
\begin{equation*}
\left\|\delta^\alpha (-)\right\|_\mathbf{t},\, \quad\text{ for all } \alpha\in \mathbb{N}^n,
\end{equation*}
where $\delta^\alpha := \delta_1^{\alpha_1}\cdots \delta_n^{\alpha_n}$, and $\alpha:=(\alpha_1,\ldots,\alpha_n)$. \\

\begin{remark}
If $\theta$ has integer entries, $\mathcal{A}_\theta^n=\mathcal{A}_0^n$ is isomorphic to the (commutative) algebra $C^\infty(\mathbb{T}^n)$ (under pointwise multiplication). \\
\end{remark}

Looking at the particular case of a noncommutative 2--torus we have the following: We fix an irrational number $\epsilon$, set 
$\theta:=\begin{pmatrix}
0 & \epsilon \\
-\epsilon &0
\end{pmatrix}$, and consider $e_1=(1,0)$ and $e_2=(0,1)$ the standard basis for $\mathbb{Z}^2$. With the cocycle $c$ defined by $\theta$ as in \eqref{def:cocycle}, the two unitary generators $U$ (regular representation), and $V$ (irrational rotation) of  the irrational rotation $C^*$--algebra $C^*(\mathbb{Z}^2,c)$ are given by \cite[Sect.\ 1]{cohenconnes}
\begin{align*}
U(f)(x):=U_{(1,0)}(f)(x)=e^{2\pi ix}f(x),\\ 
V(f)(x):=U_{(0,1)}(f)(x)=f(x+\epsilon),
\end{align*}
for $f\in L^2(S^1)$, $x\in S^1$.
With the notation above, since
\begin{equation*}
\langle e_1,\theta e_2\rangle=
\begin{pmatrix}
1 & 0
\end{pmatrix}
\theta
\begin{pmatrix}
0 \\ 1
\end{pmatrix}=
\begin{pmatrix}
1 & 0
\end{pmatrix}
\begin{pmatrix}
0 & \epsilon \\
-\epsilon &0
\end{pmatrix}
\begin{pmatrix}
0 \\ 1
\end{pmatrix}
=\epsilon,
\end{equation*}
relation \eqref{eq:product weyl elements} reads 
\begin{equation*}
UV=e^{-2\pi i\epsilon}VU,
\end{equation*}
and we can also check that the operators are unitary: $U^*=U^{-1}$, $V^*=V^{-1}$. \\
As in \eqref{eq:actiontorus}, we can introduce the dynamical system given by the action of $\mathbb{T}^2$ on $C^*(\mathbb{Z}^2,c)$ by the 2--parameter group of automorphisms $\{\alpha_s\}_{s\in\mathbb{R}^2}$ determined by
\begin{equation*}
\alpha_s(U^rV^t):=e^{\,2\pi i\langle s,(r,t)\rangle}U^rV^t,
\end{equation*}
for any $(r,t)\in\mathbb{Z}^2$ (see \cite[Sect.\ 2.1]{connestretkoff}). \\

\begin{remark}
In references such as \cite{connestretkoff} and \cite{fathizadehkhalkhalisc2t} the factor $2\pi$ does not appear in the exponential, whereas in \cite{haleepongeNCTI} and \cite{levyneirapaycha} it does. This depends on the corresponding convention used for the definition of the Fourier transform of a function. \\
\end{remark}

The subalgebra $\mathcal{A}_\theta^n$ of smooth elements of $C^*(\mathbb{Z}^2,c)$ consists of those $x\in C^*(\mathbb{Z}^2,c)$ such that the mapping
\begin{align*}
\mathbb{R}^2&\to C^*(\mathbb{Z}^2,c)\\
s&\mapsto\alpha_s(x)
\end{align*}
is smooth, i.e.\ the elements $a\in \mathcal{A}_\theta^n$, $a=\sum\limits_{(r,t)\in\mathbb{Z}^2}a_{r,t}U^rV^t$, such that the coefficients form a sequence $(a_{r,t})_{(r,t)\in\mathbb{Z}^2}\in \mathcal{S}(\mathbb{Z}^2)$, i.e.\ $\{\lvert r\rvert^k\lvert t\rvert^q\lvert a_{r,t}\rvert\}$ is bounded for all $k,q>0$. The associated derivations are given by the defining relations,
\begin{align}
\delta_1(U)&=U,\ &\delta_1(V)&=0, \label{eq:delta1}\\
\delta_2(U)&=0,\ &\delta_2(V)&=V, \label{eq:delta2}
\end{align}
which are analogues of the differential operators $-i\partial_x$, $-i\partial_y$ on the smooth functions on $\mathbb{T}^2$. \\

Reference \cite{plazas} comprises a survey of the theory of noncommutative 2--tori, about the aspects that link the differential geometry and the algebraic geometry of that space.

\section{Pseudodifferential operators on the noncommutative torus $\mathcal{A}_\theta^n$}

In Equations \eqref{eq:symbolonopenset} and \eqref{eq:def pdos Rn} of the Introduction we mentioned how pseudodifferential operators on a closed manifold can be defined and how crucial it is the local representation of the manifold in the definition of the symbols of these operators. As there is no underlying geometric space on a noncommutative space, like a noncommutative torus, such a representation by local charts is not available. There are nevertheless alternative approaches to pseudodifferential calculi on that space (\cite{connestretkoff}, \cite{levyneirapaycha}, \cite{liu}), and in this section we revisit those definitions. Reference \cite{haleepongeNCTI} offers a useful list of works using pseudodifferential operators on noncommutative tori.

\subsection{Toroidal pseudodifferential operators}\label{sect: toroidal pdos}

In \cite[Chapt.\ 4]{ruzhanskyturunen} the foundations of the theory of pseudodifferential operators on the torus in terms of Fourier coefficients and discrete operations is developed, and in \cite{levyneirapaycha} this theory is adapted to noncommutative tori. \\

\begin{definition}[Compare with {\cite[Def.\ 3.3.1]{ruzhanskyturunen}}] Let $\mathcal{A}_\theta^{\mathbb{Z}^n}$ denote the space of functions from $\mathbb{Z}^n$ to $\mathcal{A}_\theta^n$. Given $j\in \{1,\ldots,n\}$, the linear map $\Delta_{j}:\mathcal{A}_\theta^{\mathbb{Z}^n}\to \mathcal{A}_\theta^{\mathbb{Z}^n}$ defined by
  \begin{equation*}
  \Delta_j (\sigma) (k) := \sigma(k+e_j)-\sigma(k) \quad\text{ for all } k\in\mathbb{Z}^n,
  \end{equation*}
  where $\{e_j\}_{1\leq j \leq n}$ is the canonical basis of $\mathbb{R}^n$, is a \emph{forward difference operator}. If $\alpha=(\alpha_1,\ldots,\alpha_n)\in \mathbb{N}^n$, we set $\Delta^\alpha:=\Delta_1^{\alpha_1}\cdots \Delta_n^{\alpha_n}$. \\
\end{definition}

If $\sigma,\tau\in \mathcal{A}_\theta^{\mathbb{Z}^n}$, then
 \begin{equation*}
  \Delta^\alpha (\sigma\tau) = \sum_{\beta\leq \alpha}\tbinom{\alpha}{\beta} \Delta^\beta(\sigma) T_{\beta} \Delta^{\alpha-\beta} (\tau),
  \end{equation*}
where $T_\beta(\sigma)(k):=\sigma(k+\beta)$, for $k\in\mathbb{Z}^n$. \\

Consider the notation given in Section \ref{section:NCTdiff}. \\

\begin{definition}[{\cite[Def.\ 3.1]{levyneirapaycha}}]\label{def:toroidal symbol}  A function $\sigma : \mathbb{Z}^n \to \mathcal{A}_\theta^n$ is a \emph{(discrete) toroidal symbol of order} $d\in \mathbb{R}$ on $\mathcal{A}_\theta^n$ ($\sigma\in S^d_{\mathcal{A}_\theta^n}(\mathbb{Z}^n)$), if for all $\alpha, \beta\in\mathbb{N}^n$, there exists $C_{\alpha,\beta}\in \mathbb{R}$ such that
 \begin{equation*}
 \left\|\delta^\alpha\Delta^\beta \sigma (k)\right\|_\mathbf{t} \leq C_{\alpha,\beta} (1+\lvert k\rvert)^{d-\lvert\beta\rvert},\quad\text{ for all } k \in \mathbb{Z}^n. 
 \end{equation*} 
\end{definition}

\begin{example} If $j\in \{1,\ldots,n\}$, the map $k\mapsto k_j U_0$ is a symbol in $S^1_{\mathcal{A}_\theta^n}(\mathbb{Z}^n)$. \\
Any element of $\mathcal{A}_\theta^n$ can be seen as a symbol in $S^0_{\mathcal{A}_\theta^n}(\mathbb{Z}^n)$ through the injection $\mathcal{A}_\theta^n\to S^0_{\mathcal{A}_\theta^n}(\mathbb{Z}^n)$ given by $a\mapsto (k\mapsto a)$. \\
\end{example}

If $\sigma\in S_{\mathcal{A}_\theta^n}^{d_1}(\mathbb{Z}^n)$ and $\tau\in S_{\mathcal{A}_\theta^n}^{d_2}(\mathbb{Z}^n)$ then for $k\in\mathbb{Z}^n$
\begin{align}
(\sigma\circ \tau)(k) &:=  \sum_{l\in \mathbb{Z}^n} \tau_{l}(k)\,  \sigma(l+k)\,  U_l  \label{eq:starprodsymbols}\\
 &=\sum_{l\in \mathbb{Z}^n} \tau_{l}(k)\,  \sum_{r\in \mathbb{Z}^n}\sigma_r(l+k)U_rU_l  \notag \\
 &=\sum_{l\in \mathbb{Z}^n} \tau_{l}(k)\, \,  \sum_{r\in \mathbb{Z}^n}\sigma_r(l+k)c(r,l)U_{r+l}, \notag
\end{align}
gives a symbol in $S_{\mathcal{A}_\theta^n}^{d_1+d_2}(\mathbb{Z}^n)$, and $\circ : S_{\mathcal{A}_\theta^n}^{d_1}(\mathbb{Z}^n)\times S_{\mathcal{A}_\theta^n}^{d_2}(\mathbb{Z}^n) \to S_{\mathcal{A}_\theta^n}^{d_1+d_2}(\mathbb{Z}^n)$ is a bilinear map. \\

In the usual notation, given $\sigma\in S_{\mathcal{A}_\theta^n}^{d}(\mathbb{Z}^n)$ and $\sigma_j\in S_{\mathcal{A}_\theta^n}^{d-j}(\mathbb{Z}^n)$ for $j\in\mathbb{N}$, we write $\displaystyle\sigma\sim\sum_{j\in\mathbb{N}}\sigma_j$ if for all $N\in\mathbb{N}$,
\begin{equation*}
\sigma-\sum_{j=0}^N\sigma_j\in S_{\mathcal{A}_\theta^n}^{d-N-1}(\mathbb{Z}^n).
\end{equation*}

Given an extension map $e$ taking symbols on $\mathbb{Z}^n$ to symbols on $\mathbb{R}^n$ (see \eqref{extension map} below), \cite[Thm.\ 4.17]{levyneirapaycha} gives the following property for the previous composition of symbols
\begin{equation}\label{eq:expansion formula product toroidal symbols}
\sigma\circ\tau\sim\sum_{\lvert\alpha\rvert\geq0}\dfrac{1}{\alpha!}\left(\partial_\xi^\alpha e(\sigma)\right)\vert_{\mathbb{Z}^n}\overline{\delta}^\alpha\tau,
\end{equation}
where $\overline{\delta}^\alpha$ represents the derivations $\delta^\alpha$ defined in Section 3, but acting on symbols on the noncommutative torus $\mathcal{A}_\theta^n$. \\

The space of all (discrete) toroidal symbols on $\mathcal{A}_\theta^n$, $S_{\mathcal{A}_\theta^n}(\mathbb{Z}^n):=\displaystyle\bigcup_{d\in \mathbb{R}} S^d_{\mathcal{A}_\theta^n}(\mathbb{Z}^n)$ generates an $\mathbb{R}$-- algebra under the map $\circ$ defined in \eqref{eq:starprodsymbols}, and the set of \emph{smoothing} toroidal symbols given by $S^{-\infty}_{\mathcal{A}_\theta^n}(\mathbb{Z}^n):=\displaystyle\bigcap_{d\in \mathbb{R}} S^d_{\mathcal{A}_\theta^n}(\mathbb{Z}^n)$ is an ideal of this algebra. Moreover, for each $d\in\mathbb{R}$, the space $S^d_{\mathcal{A}_\theta^n}(\mathbb{Z}^n)$ is a Fr\'echet space for the seminorms
\begin{equation*}
 p_{\alpha,\beta}^{(d)}(\sigma):=\sup_{k\in \mathbb{Z}^n} (1+\lvert k\rvert)^{-d+\lvert\beta\rvert}\left\|\delta^\alpha\Delta^\beta \sigma(k)\right\|_\mathbf{t}\, \text{ for all } \alpha, \beta\in\mathbb{N}^n.
\end{equation*}

Let $\displaystyle a:=\sum_{k\in \mathbb{Z}^n} a_k \,  U_k\in \mathcal{A}_\theta^n$ and $\sigma\in S^d_{\mathcal{A}_\theta^n}(\mathbb{Z}^n)$. By {\cite[Lem.\ 3.8]{levyneirapaycha}} we have that
  \begin{equation}\label{toroidalpdo}
  \text{Op}_\theta(\sigma)(a):=\sum_{k\in \mathbb{Z}^n} a_k \, \sigma(k)\,  U_k
  \end{equation}
  is absolutely summable in $\mathcal{A}_\theta^n$. \\

\begin{definition}[{\cite[Lem.\ 3.8]{levyneirapaycha}}]\label{def:Optheta}
A \emph{toroidal pseudodifferential operator of order} $d$ is a continuous linear operator $\mathcal{A}_\theta^n\to \mathcal{A}_\theta^n$ of the form $\emph{\text{Op}}_\theta(\sigma)$ for a symbol $\sigma\in S^d_{\mathcal{A}_\theta^n}(\mathbb{Z}^n)$. \\
\end{definition}

We denote by $\Psi_\theta^d(\mathbb{T}^n):=\text{Op}_\theta(S^d_{\mathcal{A}_\theta^n}(\mathbb{Z}^n))$ the space of toroidal pseudodifferential operators of order $d$, and we equip this space with the Fr\'echet topology defined on $S^d_{\mathcal{A}_\theta^n}(\mathbb{Z}^n)$. \\

\begin{proposition}[{\cite[Prop.\ 3.11]{levyneirapaycha}}] 
Let $d\in\mathbb{R}$.
\begin{enumerate}
 \item The quantization map $\text{\emph{Op}}_\theta: S^d_{\mathcal{A}_\theta^n}(\mathbb{Z}^n) \to \Psi^d_\theta(\mathbb{T}^n)$ is a bijection.
 \item The inverse (dequantization) map $\text{\emph{Op}}_\theta^{-1}$ satisfies
 \begin{equation*}
  \text{\emph{Op}}_\theta^{-1} (A) (k)= A(U_k)\, U_{-k}\  \text{ for all } A\in \Psi^d_\theta(\mathbb{T}^n),\ k\in \mathbb{Z}^n.
 \end{equation*}
\end{enumerate}
\end{proposition}
For the second statement, compare with \cite[Thm.\ 4.1.4]{ruzhanskyturunen}. \\

The space of all toroidal pseudodifferential operators, and the space of smoothing toroidal pseudodifferential operators are respectively given by
\begin{equation*}
 \Psi_\theta(\mathbb{T}^n):=\displaystyle\bigcup_{d\in\mathbb{R}} \Psi^d_\theta(\mathbb{T}^n),\qquad\Psi^{-\infty}_\theta(\mathbb{T}^n):= \displaystyle\bigcap_{d\in\mathbb{R}} \Psi^d_\theta(\mathbb{T}^n).
\end{equation*}
Then, $\text{Op}_\theta$ extends to a bijective linear map  $\text{Op}_\theta: S_{\mathcal{A}_\theta^n}(\mathbb{Z}^n)\to\Psi_\theta(\mathbb{T}^n)$ compatible with the filtration on the space of toroidal symbols, and induces a bijective linear map $\text{Op}_\theta: S_{\mathcal{A}_\theta^n}^{-\infty}(\mathbb{Z}^n)\to \Psi^{-\infty}_\theta(\mathbb{T}^n)$. \\

If $A\in \Psi^{d_1}_\theta(\mathbb{T}^n)$ and $B\in \Psi^{d_2}_\theta(\mathbb{T}^n)$, then $AB\in \Psi^{d_1+d_2}_\theta(\mathbb{T}^n)$, more precisely, if $\sigma\in S_{\mathcal{A}_\theta^n}^{d_1}(\mathbb{Z}^n)$ and $\tau\in S_{\mathcal{A}_\theta^n}^{d_2}(\mathbb{Z}^n)$ are such that $A=\text{Op}_\theta(\sigma)$ and $B=\text{Op}_\theta(\tau)$, then, by \eqref{eq:starprodsymbols} we obtain
\begin{equation}\label{eq:composition toroidal operators}
AB=\text{Op}_\theta(\sigma)\, \text{Op}_\theta(\tau) = \text{Op}_\theta(\sigma\circ\tau). 
\end{equation}

This implies that the quantization map $\text{Op}_\theta: S_{\mathcal{A}_\theta}(\mathbb{Z}^n) \to \Psi_\theta(\mathbb{T}^n)$ is an algebraic and topological isomorphism. Moreover, the space of all toroidal pseudodifferential operators on $\mathbb{T}^n$, $\Psi_\theta(\mathbb{T}^n)$, is an $\mathbb{R}$--filtered algebra under composition of operators, and $\Psi_\theta^{-\infty}(\mathbb{T}^n)$ is the ideal of smoothing operators. \\

This pseudodifferential calculus has been used in \cite{levyneirapaycha} to give a characterization of the canonical discrete sum and the noncommutative residue on discrete toroidal symbols, and of the corresponding canonical trace and the noncommutative residue on toroidal pseudodifferential operators. By means of the canonical trace, we can derive defect formulae for regularized traces, and as a consequence, the conformal invariance of the $\zeta$--function at zero of the Laplacian on noncommutative tori holds. This calculus has also been used in \cite{azzalilevyneirapaycha} to give an interpretation of the scalar curvature on noncommutative 2--tori as an (extended) noncommutative residue.

\subsection{Connes' pseudodifferential operators}\label{sect: pdos connes}

Another pseudodifferential calculus on noncommutative tori is associated to a $C^*$--dynamical system (\cite{connescstar}, \cite{baaj1, baaj2}), analogue to the standard pseudodifferential calculus on $\mathbb{R}^n$ described in Equation \eqref{eq:def pdos Rn} of the Introduction. This calculus, also known as Connes' pseudodifferential calculus, has been used to give a proof of the Gauss--Bonnet theorem for noncommutative tori (\cite{connestretkoff}, \cite{fathizadehkhalkhaligb}, \cite{fathizadehkhalkhalisc2t}), and has been studied in great detail in \cite{haleepongeNCTI}, \cite{haleepongeNCTII} and \cite{tao}. In \cite{haleepongeNCTI} there is a very good list of references that apply these pseudodifferential techniques to the differential geometric study of noncommutative tori. \\

In \cite{leschmoscovici19} and \cite{leschmoscovici16} Connes' pseudodifferential calculus is shown to be the crucial technical tool for the explicit computation of the modular Gaussian curvature for noncommutative tori, due to the effectiveness of the variational methods. In those articles, the authors also implement the Morita equivalence by the Heisenberg bimodules for noncommutative tori, and they use the associated symbol calculus in order to compute the resolvent trace expansion, or equivalently the heat trace expansion, for the relevant Laplace--type operators. For this, they describe an extension to twisted crossed products of Connes' \cite{connescstar} and Baaj's \cite{baaj1, baaj2} pseudodifferential calculus for $C^*$--dynamical systems. \\

The $C^*$--dynamical system in our case is given by $({\mathcal{A}_\theta^n},{\mathbb{T}^n},\alpha)$ for a fixed $\theta$, where $\alpha$ denotes the action of the $n$--torus given by \eqref{eq:actiontorus} (see the end of Section \ref{section:NCTdiff} for the case of noncommutative 2--tori). \\

Consider the derivations $\delta_j$ given in \eqref{eq:deltas} by the relations $\delta_j(U_j)=U_j$, $\delta_j(U_k)=0$ for $j\not=k$. By {\cite{connescstar}} and \cite[Sect.\ 4]{connestretkoff}, similar to the definition of smooth toroidal symbols (Definition \ref{def:toroidal symbol}), we have the following \\

\begin{definition}[{\cite[Def.\ 3.1]{levyneirapaycha}}]\label{def:smooth toroidal symbol}  A smooth function $\sigma : \mathbb{R}^n \to {\mathcal{A}_\theta^n}$ is a \emph{smooth toroidal symbol of order} $d\in \mathbb{R}$ on ${\mathcal{A}_\theta^n}$ ($\sigma\in S^d_{\mathcal{A}_\theta^n}(\mathbb{R}^n)$), if for all $\alpha, \beta\in\mathbb{N}^n$, there exists $C_{\alpha,\beta}\in \mathbb{R}$ such that 
\begin{equation*}
 \left\|\delta^\alpha\partial_\xi^\beta \sigma (\xi)\right\|_\mathbf{t} \leq C_{\alpha,\beta} (1+\lvert\xi\rvert)^{d-\lvert\beta\rvert},\quad\text{ for all } \xi\in\mathbb{R}^n.
\end{equation*}
\end{definition}

Let $S_{\mathcal{A}_\theta^n}(\mathbb{R}^n):=\displaystyle\bigcup_{d\in \mathbb{R}} S^d_{\mathcal{A}_\theta^n}(\mathbb{R}^n)$ be the space of all smooth toroidal symbols on $\mathcal{A}_\theta^n$.

\begin{definition}[{\cite[Def.\ 4.1]{connestretkoff}}]\label{def:connespdo}
Let $d\in\mathbb{R}$ and let $\sigma\in S^d_{\mathcal{A}_\theta^n}(\mathbb{R}^n)$ be a smooth toroidal symbol.
The \emph{pseudodifferential operator} associated to $\sigma$ is given by the oscillatory integral
\begin{equation}\label{eq:Opquantization}
\text{\emph{Op}}(\sigma)(a):=\int_{\mathbb{R}^n}\int_{\mathbb{T}^n} e^{-2\pi i\langle s,\xi\rangle}\sigma(\xi)\alpha_s(a)\,ds\,d\xi,
\end{equation}
for any $a\in \mathcal{A}_\theta^n$. The space of pseudodifferential operators on $\mathcal{A}_\theta^n$ is given by the image of $S_{\mathcal{A}_\theta^n}(\mathbb{R}^n)$ under the map $\text{\emph{Op}}$. \\
\end{definition}

It follows from the definition of the action $\alpha$ that Equation \eqref{eq:Opquantization} corresponds to the usual definition of a pseudodifferential operator on the (commutative) manifold $\mathbb{T}^n$ (see Equation \eqref{eq:def pdos Rn}). \\

\begin{remark}
In references \cite{haleepongeNCTI} and \cite{haleepongeNCTII} there is a detailed study of these operators and a careful explanation of the meaning of this integral. In those references, the sign of the variable $s$ is different from the sign we use here. \\
\end{remark}

For $d_1,d_2\in\mathbb{R}$, let $\sigma\in S^{d_1}_{\mathcal{A}_\theta^n}(\mathbb{R}^n)$, $\tau\in S^{d_2}_{\mathcal{A}_\theta^n}(\mathbb{R}^n)$. Then 
\begin{equation}\label{eq:product symbols connes}
(\sigma\#\tau)(\xi):=\int_{\mathbb{R}^n}\int_{\mathbb{T}^n} e^{-i\langle t,\eta\rangle}\sigma(\xi+\eta)\alpha_{t}(\tau(\xi))\,dt\,d\eta, \quad \text{ for any } \xi\in\mathbb{R}^n,
\end{equation}
belongs to $S^{d_1+d_2}_{\mathcal{A}_\theta^n}(\mathbb{R}^n)$ \cite[Sect.\ 7]{haleepongeNCTII} and satisfies the following:
\begin{enumerate}
\item $\displaystyle (\sigma\#\tau)(\xi)\sim\sum_{\lvert\alpha\rvert\geq0}\dfrac{1}{\alpha!}\partial_\xi^\alpha\sigma(\xi)\delta^\alpha\tau(\xi)$.
\item $\text{Op}(\sigma)\text{Op}(\tau)\sim\text{Op}(\sigma\#\tau)$.
\end{enumerate}
Compare item 1 with \eqref{eq:expansion formula product toroidal symbols}, and item 2 with \eqref{eq:composition toroidal operators}. \\

This way, we obtain a pseudodifferential calculus on noncommutative tori analogue to the standard pseudodifferential calculus on closed manifolds. In fact, for every $d\in\mathbb{R}$, $\text{Op}(S^d_ {\mathcal{A}_0^n}(\mathbb{R}^n))$ is the standard space of pseudodifferential operators of order $d$ on the algebra of smooth functions on the commutative torus $C^\infty(\mathbb{T}^n)$. \\

One natural question is how to compare these pseudodifferential operators with the toroidal pseudodifferential operators given in the previous section. This was studied in \cite[Sect.\ 4]{levyneirapaycha}, which considers results given in \cite{ruzhanskyturunen} for periodic pseudodifferential operators; see also \cite{connestretkoff}, \cite[Prop.\ 5.9]{haleepongeNCTI}, \cite[Lem.\ 2.5]{tao}.\\

Considering the action given in \eqref{eq:actiontorus}, for $a=\displaystyle\sum_{k\in\mathbb{Z}^n}a_kU_k\in \mathcal{A}_\theta^n$, Equation \eqref{eq:Opquantization} amounts to 
\begin{equation}
\text{Op}(\sigma)(a)=\sum_{k\in\mathbb{Z}^n}a_k\sigma(k)U_k,
\end{equation}
which coincides with the definition of a toroidal pseudodifferential operator given in \eqref{toroidalpdo}.\\

In order to compare the quantization given in Definition \ref{def:Optheta} by \eqref{toroidalpdo} and the one given in Definition \ref{def:connespdo} by \eqref{eq:Opquantization}, we need to consider the extension of a symbol defined on $\mathbb{Z}^n$ to a symbol defined on $\mathbb{R}^n$. If $\sigma\in S_{\mathcal{A}_\theta^n}(\mathbb{Z}^n)$, an \emph{extension} of $\sigma$ is a symbol $\widetilde{\sigma}\in S_{\mathcal{A}_\theta^n}(\mathbb{R}^n)$ such that ${\widetilde{\sigma}}\vert_{\mathbb{Z}^n}=\sigma$, and the map 
\begin{align}
S_{\mathcal{A}_\theta^n}(\mathbb{Z}^n)/S^{-\infty}_{\mathcal{A}_\theta^n}(\mathbb{Z}^n) &\to S_{\mathcal{A}_\theta^n}(\mathbb{R}^n)/S^{-\infty}_{\mathcal{A}_\theta^n}(\mathbb{R}^n) \notag \\
\sigma&\mapsto \widetilde{\sigma} \label{extension map}
\end{align}
is a linear isomorphism (this is proved in \cite[Lem.\ 5.12]{levyneirapaycha} for complex valued symbols, but the same proof also works for $\mathcal{A}_\theta^n$--valued symbols).  Compare with \cite[Thm.\ 4.5.3]{ruzhanskyturunen}. \\

Moreover, for any $d\in\mathbb{Z}$
\begin{equation*}
\text{Op}_\theta(S^d_{\mathcal{A}_\theta^n}(\mathbb{Z}^n))=\text{Op}(S^d_{\mathcal{A}_\theta^n}(\mathbb{R}^n)),
\end{equation*}
where, on the left hand side we have the quantization given in Equation \eqref{toroidalpdo}, and on the right hand side the quantization is the one given in Equation \eqref{eq:Opquantization}  (compare with \cite[Cor.\ 4.6.13]{ruzhanskyturunen}). This way, the toroidal pseudodifferential calculus and the Connes' pseudodifferential calculus on the noncommutative torus $\mathcal{A}_\theta^n$ coincide. \\

\subsection{Deformed pseudodifferential operators}\label{sect:deformed pdos}

The deformation of manifolds along an isometric action of a real $n$--torus $\mathbb{T}^n$ with $n\geq2$ is studied in \cite{brainlandisujlekom}. That construction characterizes noncommutative toric manifolds, such as the noncommutative 4--sphere, in terms of the algebra of its smooth functions (see Remark \ref{rem:algebrafunctns}). This functorial procedure deforms not just the algebra itself but any associated $\mathbb{T}^n$--equivariant construction on the manifold, in particular its differential and metric structures, and it allows to apply the $\theta$--deformation to the whole pseudodifferential calculus on closed manifolds. \\

On the other hand, the definition of pseudodifferential operators on noncommutative tori as a deformation of the usual operators on a closed manifold can be found in \cite{liu}. In order to apply the deformation theory given in Section \ref{sect:rieffel deformation} to noncommutative tori, both the symbol map and the quantization map in the calculus have to be equivariant with respect to the torus action. This leads the author of \cite{liu} to work with a global pseudodifferential calculus on closed manifolds in which all the ingredients are given in a coordinate--free way, and he follows Widom's framework from \cite{widomfam} and \cite{widomcomp}. Now we closely follow \cite{liu}. \\

We start by introducing the deformation of operators acting on Hilbert spaces following the construction given in Section \ref{sect:deformation cstar algebra}. If $H_1$ and $H_2$ are two Hilbert spaces, $B(H_1,H_2)$ denotes the space of all bounded operators from $H_1$ to $H_2$, and if $H_1=H_2$, then $B(H_1)$ denotes the space $B(H_1,H_1)$.\\

Let $H_1$ and $H_2$ be two Hilbert spaces which are both strongly continuous unitary representations of $\mathbb{T}^n$. We denote those representations by $t\mapsto R_t \in B(H_1)$ and $t \mapsto \tilde{R}_t \in B(H_2)$ respectively, where $t \in \mathbb{T}^n$. \\ 

As in Section \ref{sect:Afrechet} any function $f\in H_1$ admits an \emph{isotypical decomposition}
\begin{equation*}
f=\sum_{l\in\mathbb{Z}^n}f_l,
\end{equation*}
where $f_l:=\displaystyle\int_{\mathbb{T}^n}R_t(f)e^{-2\pi i\langle l, t\rangle}\,dt$ belongs to the $l$--th spectral subspace (or $l$--th isotypic component), given by the torus action $R$ (see \eqref{eq:spectralsubspace}). \\

The space $B(H_1,H_2)$ becomes a $\mathbb{T}^n$--module via the adjoint action:
\begin{equation}
\label{eq:def Ad}
P\in B(H_1,H_2)\mapsto \text{Ad}_t(P):=\tilde{R}_tPR_{-t},\quad t\in\mathbb{T}^n.
\end{equation}

\begin{remark}
The action of the torus $\mathbb{T}^n$ on $A_\theta^n$ given in \eqref{eq:actiontorus} can also be seen in this way as an adjoint action on $B(L^2(\mathbb{T}^n))$ \cite[Sect.\ 2.2]{haleepongeNCTII}. \\
\end{remark}

With the notation of Section \ref{sect:Afrechet} we denote by $B(H_1,H_2)^\infty$ the space of all the smooth elements in $B(H_1,H_2)$ for the action of $\mathbb{T}^n$ given by \eqref{eq:def Ad}. Then $B(H_1,H_2)^\infty$ is a Fr\'echet space whose topology is defined by the following seminorms, which allows us to prove the continuity of the torus action 
\begin{equation*}
q_j(P):=\sum_{\lvert\beta\rvert\leq j}\dfrac{1}{\beta!}\left\|\partial_t^\beta \text{Ad}_t(P)\right\|_{B(H_1,H_2)}, \quad \text{ for all }j\in\mathbb{N}.
\end{equation*}

Any $\mathbb{T}^n$--smooth operator $P$ admits an isotypical decomposition
\begin{equation*}
P=\sum_{k\in\mathbb{Z}^n}P_k,
\end{equation*}
where $P_k:=\displaystyle\int_{\mathbb{T}^n}\text{Ad}_t(P)e^{-2\pi i\langle k, t\rangle}\,dt$, and the sequence of operator norms $\{\|P_k\|\}_{k\in\mathbb{Z}^n}$ decays faster than any polynomial in $k$. In particular, the convergence of the infinite sum is absolute with respect to the operator norm in $B(H_1,H_2)$. \\

For a fixed $n\times n$ skew symmetric matrix $\theta$, we recall the associated cocycle $c(k,l)=e^{-\pi i\langle k,\theta l\rangle}$ given in \eqref{def:cocycle}. The deformation map
\begin{equation*}
\pi^\theta: B(H_1,H_2)^\infty \to B(H_1,H_2)^\infty,
\end{equation*}
is defined by
\begin{equation}\label{eq:deformation map pitheta}
\pi^\theta(P)(f):=\sum_{k,l\in\mathbb{Z}^n}c(k,l)P_k(f_l),
\end{equation} 
with the isotypical decompositions $P=\sum\limits_{k\in\mathbb{Z}^n}P_k \in B(H_1,H_2)^\infty$ and $f=\sum\limits_{l\in\mathbb{Z}^n}f_l\in H_1$. \\

According to \cite[Lem.\ 2.9]{liu}, for any multiindex $\mu$ we can find an integer $l$ large enough such that there exists a constant $C_\mu$ with
\begin{equation*}
\left\|\partial_t^\mu(\text{Ad}_t(\pi^\theta(P)))\right\|\leq C_\mu q_l(P), \quad \text{ for all }P\in B(H_1,H_2)^\infty.
\end{equation*}

\begin{proposition}[{\cite[Lem.\ 2.10, Prop.\ 2.12]{liu}}]\label{prop: deformation map operators}
Let $H$ be a Hilbert space. The deformation map
\begin{equation*}
\pi^\theta:(B(H)^\infty,\times_\theta)\to(\pi^\theta(B(H)^\infty),\circ)
\end{equation*}
is a $*$--algebra isomorphism, namely, for any $P,Q\in B(H)^\infty$,
\begin{equation}\label{eq: pitheta composition}
\pi^\theta(P\times_\theta Q)=\pi^\theta(P)\circ\pi^\theta(Q),\qquad \pi^\theta(P^*)=\pi^\theta(P)^*
\end{equation}
where,
\begin{equation}\label{eq:prodthetaoperators}
P\times_\theta Q:=\displaystyle\sum\limits_{k,l\in\mathbb{Z}^n}c(k,l)P_k\circ Q_l=\displaystyle\sum\limits_{k,l\in\mathbb{Z}^n}e^{-\pi i\langle k,\theta l\rangle}P_k\circ Q_l.
\end{equation}
\end{proposition}
Here we can compare \eqref{eq:prodthetaoperators} with \eqref{eq:product of eigenvectors} for $J=\theta/2$ as before. \\

Let $M$ be a closed manifold which admits an $n$--torus action: $(t,x)\mapsto t\cdot x$ for $t\in\mathbb{T}^n$ and $x\in M$ (in \cite[Sect.\ 2.3]{liu} it is considered the action of the inverse of $t$ on $x$). We obtain an $n$--torus action on $C^\infty(M)$:
\begin{equation}\label{eq:torusactionCinftyM}
\alpha_t(f)(x):=f(t\cdot x),\qquad f\in C^\infty(M).
\end{equation}
This action defines a $*$--automorphism of $C^\infty(M)$, so that $\mathbb{T}^n$ acts smoothly on $C^\infty(M)$ with respect to the smooth Fr\'echet topology given by the seminorms
\begin{equation*}
\lvert f\rvert_{l,K}:=\sup_{\lvert\alpha\rvert\leq l,\ x\in K}\left\lvert\partial^\alpha f(x)\right\rvert,
\end{equation*}
for $l\in\mathbb{N}$ and $K\subset M$ compact. The {pointwise multiplication} is jointly continuous and therefore we can deform it to $\times_\theta$ as in Definition \ref{def:deformedproduct} with respect to a skew symmetric matrix $\theta$. The new algebra 
\begin{equation*}
C^\infty(M_\theta):=(C^\infty(M),\times_\theta)
\end{equation*}
plays the role of smooth coordinate functions on a \emph{noncommutative manifold} $M_\theta$. \\

\begin{example}
In Section \ref{section:NCTalg} we identified the noncommutative torus $\mathcal{A}_\theta^n$ with the space obtained after a deformation of the product on $C^\infty(\mathbb{T}^n)$, so that $\mathcal{A}_\theta^n=C^\infty(\mathbb{T}^n_\theta)$.  Considering the functions $E_k$ on the torus $\mathbb{T}^n$, with $k\in\mathbb{Z}^n$, and the action $\alpha$ described in \eqref{eq:action torus Ek}, the representation given in \eqref{eq:decomposition f torus} 
\begin{equation*}
\sum\limits_{k\in\mathbb{Z}^n}f_{k}E_k
\end{equation*}
becomes the isotypical decomposition of $f\in C^\infty(\mathbb{T}^n_\theta)$. Moreover, the product $\times_\theta$ is the one given in \eqref{eq:product times theta}. \\
\end{example}

It is possible to lift the torus action to $T^*M$ by the natural extension of diffeomorphisms $\varphi\mapsto\varphi^*$, where $\varphi^*$ is the differential of $\varphi$ (see \cite[p.\ 525]{liu}). Thus, the cotangent bundle of the noncommutative manifold $M_\theta$ is given by the deformed algebra 
\begin{equation*}
C^\infty(T^*M_\theta):=(C^\infty(T^*M),\times_\theta).
\end{equation*}

Now, remember the definition of a symbol on an open subset of $\mathbb{R}^n$ given in \eqref{eq:symbolonopenset} and compare it with \cite[Sect.\ 3]{widomcomp} and with \cite[Eq.\ (2.29)]{liu}: \\
\begin{definition}[See e.g.\ {\cite[Def.\ 1.1]{shubin}}]
A smooth function $\sigma:\mathbb{R}^n\times\mathbb{R}^n\to\mathbb{C}$ is a \emph{symbol of order} $d\in\mathbb{R}$ on $\mathbb{R}^n$ ($\sigma\in S^d(\mathbb{R}^n)$), if for all $\alpha,\beta\in\mathbb{N}^n$, there exists $C_{\alpha,\beta}\in\mathbb{R}$ such that
\begin{equation*}
\left\lvert\partial_x^\alpha\partial_\xi^\beta \sigma(x,\xi)\right\rvert\leq C_{\alpha,\beta}(1+\lvert\xi\rvert)^{d-\lvert\beta\rvert},\quad \text{ for all }\xi\in\mathbb{R}^n,
\end{equation*}
uniformly in compact $x$--sets. For a manifold $M$ the space $S^d(M)$ consists of functions satisfying this in terms of local coordinates. We also use the notation $S(M)=\bigcup\limits_{d\in\mathbb{R}}S^d(M)$ and $S^{-\infty}(M)=\bigcap\limits_{d\in\mathbb{R}}S^d(M)$. \\
\end{definition}

\begin{remark} 
The notation $S(M)$ will make sense below with Definition \ref{def:Symp and OP Liu}, when invariantly defined symbols (in a coordinate--free way) introduced by H.\ Widom in \cite{widomcomp} are considered. \\
\end{remark}

With the notation of Equation \eqref{eq:torusactionCinftyM}, for any $t\in\mathbb{T}^n$ let $\alpha_t$ denote the torus action on $C^\infty(T^*M)$. Then, for any $t\in\mathbb{T}^n$, and for any symbol $\sigma\in S^d(M)$ of order $d$, $t\mapsto \alpha_t(\sigma)$ is a $C^\infty(T^*M)$ --valued function. In local coordinates
\begin{equation*}
\partial_t^\gamma \alpha_t(\sigma)=\alpha_t\left(\sum_j\partial_x^{\alpha_j}\partial_{\xi}^{\beta_j}\sigma\right),
\end{equation*}
where $\gamma, \alpha_j, \beta_j$ are multiindices. Thus, $\alpha_t(\sigma)$ belongs to $S^d(M)$ and the torus action is smooth.  \\ 

Therefore, on the filtered algebra $S(M)$, we can deform {the pointwise multiplication}
\begin{equation*}
m:S^{d_1}(M)\times S^{d_2}(M)\to S^{d_1+d_2}(M),
\end{equation*}
by
\begin{align*}
m_\theta:S^{d_1}(M)\times S^{d_2}(M)&\to S^{d_1+d_2}(M)\\
(\sigma,\tau)&\mapsto \sum_{k,l\in\mathbb{Z}^n}c(k,l)m(\sigma_k,\tau_l),
\end{align*}
where $\displaystyle\sigma=\sum_{k\in\mathbb{Z}^n}\sigma_k$, $\displaystyle\tau=\sum_{l\in\mathbb{Z}^n}\tau_l$ are the isotypical decompositions of $\sigma$ and $\tau$ respectively. So, for $q\in T^*M$
\begin{align}
m_\theta(\sigma,\tau)(q)&:= \sum_{k,l\in\mathbb{Z}^n}c(k,l)m(\sigma_k,\tau_l)(q) \label{eq:product mtheta} \\ 
&=\sum_{l\in\mathbb{Z}^n}\tau_l(q)\sum_{k\in\mathbb{Z}^n}c(k,l)\sigma_k(q). \notag
\end{align}
(compare with Equation \eqref{eq:starprodsymbols}). Thus, the deformed version of the space of symbols on $M$ is denoted by
\begin{equation}\label{eq:symbols on Mtheta}
S(M_\theta):=(S(M),m_\theta).
\end{equation}

\begin{remark}
Looking at the action of the torus given by $\alpha_v(U_k):=e^{2\pi i\langle v, k\rangle}U_k$, the isotypical decomposition of a symbol  $\tau\in S_{\mathcal{A}_\theta^n}(\mathbb{Z}^n)$ corresponds to the expression $\displaystyle\tau(k)=\sum_{l\in\mathbb{Z}^n}\tau_l(k)U_l$ for any $k\in\mathbb{Z}^n$. \\
\end{remark}

\begin{remark}
The product that is deformed here is the pointwise multiplication and not the star product on symbols that corresponds to the symbol of the composition of the associated operators (see below in Proposition \ref{prop:product of deformed operators and symbols} the definition of the star product on symbols $\bigstar$ that is deformed to get $\bigstar_\theta$). \\
\end{remark}

As before, let $M$ be a closed manifold with a torus action which gives a diffeomorphism on $M$. A pseudodifferential operator on $M$ is an operator from $C^\infty(M)$ to $C^\infty(M)$ which locally is of the form \eqref{eq:def pdos Rn}. The author of \cite{liu} applies the deformation of operators described above to $\Psi(M)$, the algebra of standard pseudodifferential operators acting on $C^\infty(M)$. Thus, the algebra of pseudodifferential operators on the noncommutative manifold $M_\theta$ coincides with the image of $\Psi(M)$ under the deformation map $\pi^\theta$ defined in \eqref{eq:deformation map pitheta}. \\

Let $\Psi^d(M)$ be the space of pseudodifferential operators of order $d$, that is, $\Psi^d(M)$ is the space of pseudodifferential operators whose symbol, when localized on some open chart, belongs to $S^d(M)$. Let $\{\mathcal{H}_s\}_{s\in\mathbb{R}}$ be the associated Sobolev spaces of the manifold $M$. A pseudodifferential operator $P\in\Psi^d(M)$ is a continuous linear operator $P:C^\infty(M)\to C^\infty(M)$ which admits a bounded extension from $\mathcal{H}_d$ to $\mathcal{H}_0=L^2(M)$. \\

\begin{proposition}[See {\cite[Prop.\ 4.1 and Cor.\ 4.2]{liu}}]
Let $P$ be an operator, $d\in\mathbb{Z}$, and let $\{F_1,\ldots,F_l\}$ be an arbitrary finite collection of first order differential operators. Then $P$ belongs to the space $\Psi^d(M)$ if and only if $[F_l,\cdots,[F_1,P]\cdots]$ is a bounded operator from $\mathcal{H}_d$ to $\mathcal{H}_0$. \\
\end{proposition}

Let $j\in\mathbb{Z}$ and let $\mathcal{F}:=\{F_1,\ldots,F_l\}$ be a finite collection of first order differential operators. Define the following seminorms on $\Psi^d(M)$:
\begin{equation*}
\|P\|_{(j,\mathcal{F})}:=\|[F_l,\cdots,[F_1,P]\cdots]\|_{j+d,j},
\end{equation*}
where $\|\ \|_{j+d,j}$ is the operator norm from $\mathcal{H}_{j+d}(M)\to\mathcal{H}_j(M)$. For all $d\in\mathbb{Z}$, the seminorms $\|\ \|_{(j,\mathcal{F})}$ make $\Psi^d(M)$ into a Fr\'echet space. With the notation of Equation \eqref{eq:def Ad}, the function
\begin{align*}
\mathbb{T}^n&\to\Psi^d(M)\\
t&\mapsto \text{Ad}_t(-)
\end{align*}
is smooth with respect to the Fr\'echet topology of $\Psi^d(M)$ (see \cite[Prop.\ 4.4]{liu}). Moreover, from the smoothness of this torus action, the right hand side of the isotypical decomposition
\begin{equation*}
P=\sum_{k\in\mathbb{Z}^n}P_k, \quad \text{ for }P\in\Psi^d(M),
\end{equation*}
converges to $P$ with respect to such Fr\'echet topology. Additionally, for each seminorm $\|\ \|_{(j,\mathcal{F})}$, the sequence $\{\|P_k\|_{(j,\mathcal{F})}\}_{k\in\mathbb{Z}^n}$ is of rapid decay in $k$. \\

As in Equation \eqref{eq:deformation map pitheta}, fixing an $n\times n$ skew symmetric matrix $\theta$, the definition of the deformation map $\pi^\theta$ is extended to pseudodifferential operators of all orders, namely, for all $P\in\Psi(M)=\bigcup\limits_{d\in\mathbb{Z}}\Psi^d(M)$, $\pi^\theta(P)$ is given by
\begin{equation*}
\pi^\theta(P)(f)=\sum_{k,l\in\mathbb{Z}^n}c(k,l)P_k(f_l), \quad \text{ for all }f=\sum_{l\in\mathbb{Z}^n}f_l\in C^\infty(M).
\end{equation*}  
For each $d\in\mathbb{Z}$, we set $\Psi^d(M_\theta):=\pi^\theta(\Psi^d(M))$ and $\Psi(M_\theta):=\pi^\theta(\Psi(M))$.\\

The filtered algebra $\Psi(M)$ admits the following deformation: As in \eqref{eq:prodthetaoperators}, for any $P\in\Psi^{d_1}(M)$ and $Q\in\Psi^{d_2}(M)$, the $\times_\theta$ multiplication is well defined
\begin{align*}
\times_\theta:\Psi^{d_1}(M)\times\Psi^{d_2}(M)&\to\Psi^{d_1+d_2}(M)\\
(P,Q)&\mapsto P\times_\theta Q=\sum_{k,l\in\mathbb{Z}^n}c(k,l)P_k\circ Q_l.
\end{align*} 
The $\times_\theta$ multiplication is compatible with the original $*$--operation in $\Psi(M)$, namely
\begin{equation*}
(P\times_\theta Q)^*=Q^*\times_\theta P^*,
\end{equation*} 
and therefore $(\Psi(M),\times_\theta)$ is a filtered $*$--algebra. Thus, the deformation map $\pi^\theta$ makes $(\Psi(M_\theta),\circ)$ into a filtered $*$--algebra, and it satisfies the equations in \eqref{eq: pitheta composition} from Proposition \ref{prop: deformation map operators}, so that the deformation map $\pi^\theta$ is a filtered $*$--algebra isomorphism between $(\Psi(M),\times_\theta)$ and $(\Psi(M_\theta),\circ)$. \\

According to \cite[Sect.\ 4]{liu}, the space of smoothing operators is invariant under the deformation map, i.e.
\begin{equation*}
\pi^\theta(\Psi^{-\infty}(M))=\Psi^{-\infty}(M).
\end{equation*}
Therefore, the deformation of the algebra of pseudodifferential operators modulo the smoothing operators is the quotient
\begin{align*}
\Psi(M_\theta)/\Psi^{-\infty}(M)&=\pi^\theta(\Psi(M))/\Psi^{-\infty}(M)\\
&=\pi^\theta\big(\Psi(M)/\Psi^{-\infty}(M)\big).
\end{align*}

In the late seventies, H.\ Widom gave the description of a complete symbolic calculus for pseudodifferential operators in the scalar case \cite{widomfam}, and on sections of vector bundles over compact manifolds \cite{widomcomp}. In \cite[Thm.\ 2.7]{getzler}, that calculus is used in order to derive an explicit formula for composition of principal symbols, which permits a complete calculation of the index of the Dirac operator on a spin manifold. \\

Let $M$ be a closed manifold endowed with a torsion free connection $\nabla$. We will assume that the connection is invariant with respect to the torus action. The definition of Widom's complete symbolic calculus involves the tensor calculus on $M$, that consists of the pointwise tensor product, the contraction between tensor fields and the connection. If the manifold $M$ admits a torus action, it is possible to obtain a deformation of that tensor calculus to use the above results (see \cite[Sect.\ 3]{liu}). For example, any differential operator $P$ on $M$ can be defined in a coordinate--free way in the form of a finite sum
\begin{equation*}
P(f):=\sum_\alpha\rho_\alpha\cdot\left(\nabla^\alpha f\right),
\end{equation*}
where each $\rho_\alpha$ is a contravariant tensor field such that the contraction $\rho_\alpha\cdot\left(\nabla^\alpha f\right)$ gives rise to a smooth function on $M$. If we replace the contraction $\cdot$ by the deformed version $\cdot_\theta$, we get back the deformed operator $\pi^\theta(P)$. \\

The function $l(x,\xi,y):=\langle x-y,\xi\rangle$ in \eqref{eq:def pdos Rn} can be generalized to manifolds. In fact, by \cite[Prop.\ 2.1]{widomcomp}, there exists a function $\ell\in C^\infty(T^*M\times M)$ such that $\ell(\cdot,x)$ is, for each $x\in M$, linear on the fibers of $T^*M$, and such that for each $v\in T^*M$
\begin{equation*}
\partial^k\ell(v,x)\vert_{x=\pi(v)}=\begin{cases}
v, & \text{ if } k=1, \\
0, & \text{ if } k\not=1,
\end{cases}
\end{equation*}
where $\pi:T^*M\to M$ is the projection of the cotangent bundle of $M$, and $\partial^k$ represents the symmetrized $k$--th covariant derivative. \\

\begin{definition}[{\cite[Def.\ 5.1, 5.2]{liu}}]\label{def:Symp and OP Liu}
Let $M$ and $\ell$ be as above. Denote by $\psi\in C^\infty(M\times M)$ a cut--off function such that $\psi\equiv1$ on a small neighborhood of the diagonal, and such that $\psi(x,y)\not=0$ implies that $d_y\ell(\xi_x,y)\not=0$ for all $\xi_x\not=0$.
\begin{enumerate}
\item For any pseudodifferential operator $P:C^\infty(M)\to C^\infty(M)$, the symbol $\emph{\text{Symb}}(P)$ of $P$ is a smooth function on $T^*M$:
\begin{equation*}
\emph{\text{Symb}}(P)(\xi_x)=P\psi(x,y)e^{i\ell(\xi_x,y)}\vert_{y=x},
\end{equation*}
where the operator $P$ acts on the $y$ variable.
\item For any $f\in C^\infty(M)$, set $y:=\exp_xY$. The quantization map $\emph{\text{OP}}:S^d(M)\to\Psi^d(M)$, $d\in\mathbb{Z}$, is defined by
\begin{equation*}
\emph{\text{OP}}(\sigma)(f)(x):=\int_{T_x^*M}\int_{T_xM}e^{-i\langle \xi_x,Y\rangle}\sigma(\xi_x)\psi(x,y)f(\exp_xY)\,dY\,d\xi_x,
\end{equation*}
where $m=\dim(M)$, $\langle \xi_x,Y\rangle$ is the canonical pairing: $T_x^*M\times T_xM\to\mathbb{R}$, and $dY$, $d\xi_x$ denote the normalized densities that are dual to each other.
\end{enumerate}
\end{definition}

Furthermore, in \cite[Sect.\ 3]{widomcomp}, it is proved that for any $P\in \Psi^d(M)$, the symbol $\text{Symb}(P)$ belongs to the space $S^d(M)$. The choice of different $\psi$ and $\ell$ gives rise to a symbol differing from $\text{Symb}(P)$ in a smoothing symbol. Moreover, the maps $\text{Symb}$ and $\text{OP}$ are inverse to each other up to smoothing operators, resp.\ smoothing symbols. In fact, let us consider the usual notation: If $P,Q\in\Psi(M)$ are operators such that $P-Q\in\Psi^{-\infty}(M)$, we write $P\sim Q$, and similarly, if $\sigma,\tau\in S(M)$ are symbols such that $\sigma-\tau\in S^{-\infty}(M)$, we write $\sigma\sim\tau$. Following \cite[Sect.\ 5.3]{liu}, the symbol calculus depends on the choice of a linear connection and therefore, assuming that the torus acts as an affine transformation so that the connection is preserved, both the symbol map and the quantization map are $\mathbb{T}^n$--equivariant. The Fr\'echet topologies on symbols and pseudodifferential operators are constructed in such a way that both the symbol map and the quantization map are continuous. As a consequence, for any pseudodifferential operator $P=\sum_{k\in\mathbb{Z}^n}P_k$ and its symbol $\sigma=\sum_{k\in\mathbb{Z}^n}\sigma_k$ with their isotypical decompositions, we have
\begin{enumerate}
\item $\text{Symb}(P_k)\sim\text{Symb}(P)_k\sim \sigma_k$.
\item $\text{OP}(\sigma_k)\sim\text{OP}(\sigma)_k\sim P_k$.
\item $\text{Symb}(P)\sim\sum_{k\in\mathbb{Z}^n}\sigma_k$.
\item $\text{OP}(\sigma)\sim\sum_{k\in\mathbb{Z}^n}P_k$. \\
\end{enumerate}

In \cite[Sect.\ V]{widomfam} and \cite[Prop.\ 3.6]{widomcomp} H.\ Widom gave the analogue of the symbol of the composition of two operators on manifolds with a given connection. Given two pseudodifferential operators $P$ and $Q$, with symbols $\sigma$ and $\tau$ respectively, the symbol of the composition $P\circ Q$ is given by an asymptotic product
\begin{equation}\label{eq:bigstar product symbols}
\sigma\bigstar\tau\sim\sum_{j=0}^\infty a_j(\sigma,\tau),
\end{equation}
where $a_j(\ ,\ )$ are bidifferential operators reducing the total degree by $j$, which are given in terms of tensor products and contraction between tensor fields. \\

\begin{proposition}[{\cite[Prop.\ 5.7]{liu}}]\label{prop:product of deformed operators and symbols}
Let $M$ be a closed manifold with an $n$--torus action so that the previous deformation machinery applies. Let $P=\sum_{k\in\mathbb{Z}^n}P_k,\ Q=\sum_{l\in\mathbb{Z}^n}Q_l$ be pseudodifferential operators with symbols $\sigma=\sum_{k\in\mathbb{Z}^n}\sigma_k,\ \tau=\sum_{l\in\mathbb{Z}^n}\tau_l$ respectively, and let $\pi^\theta(P)$ and $\pi^\theta(Q)$ be the corresponding deformed pseudodifferential operators. Then
\begin{equation}\label{eq:pitheta OP composition symbols}
\pi^\theta(P)\circ\pi^\theta(Q)\sim\pi^\theta(\text{\emph{OP}}(\sigma\bigstar_\theta \tau)),
\end{equation}
with
\begin{equation}
\sigma\bigstar_\theta \tau\sim \sum_{j=0}^\infty a_j(\sigma, \tau)_\theta,
\end{equation}
where
\begin{equation}
a_j(\sigma, \tau)_\theta:=\sum_{k,l\in\mathbb{Z}^n}c(k,l)a_j(\sigma_k, \tau_l),
\end{equation}
and $a_j(\ ,\ )$ are the bidifferential operators appearing in the star product $\bigstar$ of symbols given in \eqref{eq:bigstar product symbols}. \\
\end{proposition}

The proof of this follows from Equation \eqref{eq: pitheta composition}, which states that $\pi^\theta(P)\circ\pi^\theta(Q)=\pi^\theta(P\times_\theta Q)$. Thus, since $\pi^\theta$ is an isomorphism, $\sigma\bigstar_\theta \tau$ should be $\text{Symb}(P\times_\theta Q)$ up to a smoothing symbol. Formally
\begin{align*}
\text{Symb}(P\times_\theta Q)&\sim\sum_{k,l\in\mathbb{Z}^n}c(k,l)\text{Symb}(P_k\circ Q_l)\\
&\sim \sum_{k,l\in\mathbb{Z}^n}c(k,l)\sigma_k \bigstar \tau_l\\
&\sim \sum_{k,l\in\mathbb{Z}^n}\sum_{j=0}^\infty c(k,l)a_j(\sigma_k, \tau_l)\\
&\sim \sum_{j=0}^\infty \sum_{k,l\in\mathbb{Z}^n}c(k,l)a_j(\sigma_k, \tau_l)\\
&\sim \sum_{j=0}^\infty a_j(\sigma, \tau)_\theta.
\end{align*}
In \cite[Eq.\ (5.20)]{liu} there is an explicit description of the terms $a_j(\ ,\ )_\theta$. \\

\begin{remark}
As we mentioned before, the product $\bigstar$ between two symbols is the symbol of the composition of the associated operators up to a smoothing symbol. Moreover, given two symbols $\sigma,\tau\in S(M)$,
\begin{equation}
\sigma\bigstar_\theta\tau-m_\theta(\sigma,\tau)\sim\sum_{j=1}^\infty a_j(\sigma, \tau)_\theta.
\end{equation}
We recall from Equation \eqref{eq:symbols on Mtheta} that $S(M_\theta):=(S(M),m_\theta)$ is defined by using the deformation of the pointwise multiplication and not the deformation of the product $\bigstar$. \\
\end{remark}

We end up this section by mentioning that this deformed pseudodifferential calculus has been used in \cite{liu} to derive a general expression for the modular curvature with respect to a conformal change of metric on even dimensional toric noncommutative manifolds. \\

\section{Noncommutative tori as the underlying noncommutative algebras of certain spectral triples}

In noncommutative geometry, a geometric space is usually described in terms of a spectral triple (see \cite{careypr}, \cite{connesmoscovici}). For a study of this concept in physics we refer the reader to \cite{suijlekom}. In this section we recall this definition and study the example where a noncommutative torus is the corresponding underlying noncommutative algebra. \\

\begin{definition}[{\cite[Def.\ 2.19]{careypr}}]
A \emph{spectral triple} $(\mathcal{A}, \mathcal{H}, D)$ consists of
\begin{enumerate}
\item an involutive algebra $\mathcal{A}$,
\item a complex Hilbert space $\mathcal{H}$,
\item an unbounded self--adjoint operator $D:\text{Dom}(D)\subseteq\mathcal{H}\to\mathcal{H}$, 
\end{enumerate}
such that 
\begin{enumerate}
\item[a.] $\mathcal{A}$ is represented on $\mathcal{H}$: there exists a $*$--homomorphism
\begin{equation*}
\pi:\mathcal{A}\to B(\mathcal{H}).
\end{equation*}
\item[b.] For all $a\in\mathcal{A}$, $\pi(a)\text{Dom}(D)\subseteq\text{Dom}(D)$.
\item[c.] For all $a\in\mathcal{A}$, the densely defined operator $[D,\pi(a)]:=D\pi(a)-\pi(a)D$ extends to a bounded operator on $\mathcal{H}$.
\item[d.] For all $a\in\mathcal{A}$, $\pi(a)(1+D^2)^{-1/2}$ is a compact operator. 
\end{enumerate}
\end{definition}

\begin{example}\label{ex:NC2T}
As a noncommutative space, a noncommutative torus can be seen as the underlying noncommutative algebra of certain spectral triples, and here we recall two of those constructions: 
\begin{enumerate}
\item Consider the noncommutative 2--torus $\mathcal{A}_\theta^2$ described in Section \ref{section:NCTdiff}:
\begin{equation*}
\mathcal{A}_\theta^2=\left\{\sum_{r,t\in \mathbb{Z}} a_{r,t} U^rV^t:\  (a_{r,t})_{r,t}\in \mathcal{S}(\mathbb{Z}^2)\right\},
\end{equation*}
together with the map $\mathbf{t}$ given in \eqref{eq:trace t} and the inner product given in \eqref{eq:inner product t} with associated norm $\|\ \|_\mathbf{t}$. Equipped with that inner product, $\mathcal{A}_\theta^2$ becomes a pre--Hilbert space. Let $L^2(\mathcal{A}_\theta^2,\mathbf{t})$ be the corresponding completion. We set 
\begin{equation}\label{eq:Hilbert space Atheta}
\mathcal{H}:=L^2(\mathcal{A}_\theta^2,\mathbf{t})\oplus L^2(\mathcal{A}_\theta^2,\mathbf{t}).
\end{equation}
Using the derivations $\delta_1, \delta_2$ given in \eqref{eq:delta1} and \eqref{eq:delta2}, we define
\begin{equation*}
D:=\begin{pmatrix}
0 & \delta_1+i\delta_2 \\
-\delta_1+i\delta_2 &0
\end{pmatrix}.
\end{equation*}
Thus, $D$ is a Dirac operator, and $(\mathcal{A}_\theta^2,\mathcal{H},D)$ is a spectral triple with the noncommutative 2--torus $\mathcal{A}_\theta^2$ as the underlying noncommutative algebra \cite[Sect.\ 3]{careypr}. \\

\item In \cite[Sect.\ 8]{savinschrohe} the authors describe a local index formula for noncommutative tori in spectral triples associated with the affine metaplectic group. We describe here their construction for the case of a noncommutative 2--torus:\\
Let $v_1,v_2\in\mathbb{C}^2$ be two vectors linearly independent over $\mathbb{Q}$; they generate the lattice
\begin{equation*}
\{n_1v_1+n_2v_2\ :\ n_1,n_2\in\mathbb{Z}\}\subseteq\mathbb{C}^2,
\end{equation*}
which is isomorphic to $\mathbb{Z}^2$. Consider the Heisenberg--Weyl operators:
\begin{equation*}
T_{v_j}u(x):=e^{i(k_jx-a_jk_j/2)}u(x-a_j),
\end{equation*}
where $a_j:=\emph{\text{Re}}(v_j)$, $k_j:=-\emph{\text{Im}}(v_j)$, for $j=1,2$. These operators are unitary and they act on the graded Hilbert space $\mathcal{H}:=L^2(\mathbb{R}^2,\Lambda(\mathbb{C}^2))$, where the identification $\Lambda(\mathbb{C}^2)\cong\Lambda(\mathbb{R}^2)\otimes\mathbb{C}$ is used. They satisfy the commutation relations
\begin{equation*}
T_{v_k}T_{v_l}=e^{-i{\emph{\text{Im}}}(v_k\overline{v_l})}T_{v_l}T_{v_k},
\end{equation*}
for $k,l=1,2$. These operators are also extended to the space of forms by the trivial action on the differentials, and that extension is denoted by the same symbol. Define the algebra of functions on a 2--dimensional noncommutative torus by
\begin{equation*}
\mathcal{A}:=\left\{\sum_\ell c_\ell T_{v_1}^{\ell_1}T_{v_2}^{\ell_2}\ :\ c_\ell\in\mathbb{C},\ \ell=(\ell_1,\ell_2)\in\mathbb{Z}^2\right\},
\end{equation*}
and consider the operator $D$ on the Schwartz space of complex valued differential forms
\begin{equation*}
D:\mathcal{S}(\mathbb{R}^2,\Lambda^{\rm{ev}}(\mathbb{R}^2)\otimes\mathbb{C})\to\mathcal{S}(\mathbb{R}^2,\Lambda^{\rm{odd}}(\mathbb{R}^2)\otimes\mathbb{C}),
\end{equation*}
defined by
\begin{equation*}
D:=\dfrac{1}{\sqrt{2}}\left(d+d^*+x\,dx\wedge+(x\,dx\wedge)^*\right),
\end{equation*}
where $d$ represents the exterior differential, $x\,dx\wedge$ is the exterior multiplication by $x\,dx=\dfrac{dr^2}{2}=\sum\limits_jx_j\,dx_j$, where $r=\lvert x\rvert$, and $d^*$ and $(x\,dx\wedge)^*$ stand for the adjoint operators. Then $(\mathcal{A},\mathcal{H},D)$ is a spectral triple with a noncommutative 2--torus as the underlying noncommutative algebra. \\
\end{enumerate}

\end{example}

\begin{remark} In references such as \cite{carotenutodabrowski}, \cite[Sect.\ 5]{essouabriiochumlevysitarz},  \cite{graciabondia}, \cite{venselaar}, we can also find descriptions of spectral triples with a noncommutative 2--torus as the underlying noncommutative algebra, used to study spin structures over the torus. For that purpose, the authors consider the same algebra $\mathcal{A}_\theta^2$ given above, but for the Hilbert space they consider a complexification of \eqref{eq:Hilbert space Atheta}. We also refer to \cite{connesmoscovicimodularcurvature} for the study of modular spectral triples on noncommutative 2--tori describing the curved geometry of these spaces.
\end{remark}

\section{Abstract pseudodifferential operators}\label{sect:abstractpdos}

In \cite{connesmoscovici}, A.\ Connes and H.\ Moscovici gave a setup for an abstract pseudodifferential calculus in order to prove a local index formula for arbitrary spectral triples satisfying special conditions (finite summability and simple dimension spectrum). In this section we recall the construction of that calculus following \cite{connesmoscovici}, \cite{higson}, \cite{higson06} and \cite{paycha}. \\

Let $(\mathcal{A},\mathcal{H},D)$ be a spectral triple. We set $\Delta:=D^2$ and for any complex number $z$
\begin{equation*}
\Delta^z:=\int_0^\infty\lambda^z\,dP_\lambda,
\end{equation*}
where $P_\lambda:=1_{(-\infty,\lambda]}(\Delta)$ stands for the spectral projection of the self--adjoint operator $\Delta$ corresponding to the interval $(-\infty,\lambda]$. The space $\mathcal{H}^z$ is obtained as the completion of $\mathcal{H}^\infty:=\bigcap\limits_{k=1}^\infty\text{Dom}(D^k)$ for the norm
\begin{equation*}
\|u\|_z:=\|\Delta^{z/2}u\|.
\end{equation*}
If $\text{Re}(z)\geq0$, then $\mathcal{H}^z$ lies in $\mathcal{H}$ and coincides with the domain $\text{Dom}(\Delta^{z/2})$. If $\text{Re}(z)<0$, then $\Delta^{z/2}$ is bounded on $\mathcal{H}$, whereas $\mathcal{H}\subsetneq\mathcal{H}^z$. For each $s\in\mathbb{R}$, we have that $\mathcal{H}^s=\text{Dom}(\lvert D \rvert^s)$ and
\begin{equation*}
\mathcal{H}^\infty=\bigcap_{s\geq0}\mathcal{H}^s,\qquad \mathcal{H}^{-\infty}:=\text{ dual of }\mathcal{H}^\infty.
\end{equation*}
Thus, we obtain a scale of Hilbert spaces (see \cite[Def.\ 10.11]{graciabondia}). \\

Let us set $\lvert D \rvert:=\Delta^{1/2}$, and denote by $\delta$ the unbounded derivation of the space $B(\mathcal{H})$ given by the commutator with $\lvert D \rvert$:
\begin{equation*}
\delta(T):=[\lvert D \rvert,T], \quad \text{ for all }T\in B(\mathcal{H}).
\end{equation*}

Let us assume the following smoothness condition on $\mathcal{A}$, which implies that the space $\mathcal{H}^\infty$ is stable under left multiplication by elements of $\mathcal{A}$ (\cite[Rem.\ 3.24]{higson}): \\
\begin{definition}[{\cite[Def.\ 10.10]{graciabondia}, \cite[Def.\ 4.25]{higson06}}]
The spectral triple $(\mathcal{A},\mathcal{H},D)$ is \emph{regular} if for all $a\in\mathcal{A}$, both $a$ and $[D,a]$ are in the domains of all powers of the derivation $\delta$: 
\begin{equation*}
a,[D,a]\in\bigcap\limits_{n\in\mathbb{N}}\text{\emph{Dom}}(\delta^n).
\end{equation*}
\end{definition}

Let us define 
\begin{equation*}
\text{OP}^0:=\bigcap\limits_{n\in\mathbb{N}}\text{Dom}(\delta^n),
\end{equation*}
and for any real number $s$, let $\text{OP}^s$ be the set of operators defined by:
\begin{equation*}
P\in\text{OP}^s \text{ if and only if }\lvert D \rvert^{-s}P\in\text{OP}^0.
\end{equation*}
In other words
\begin{equation*}
\text{OP}^s:=\lvert D \rvert^{s}\text{OP}^0.
\end{equation*}
Thus, for all $s\in\mathbb{R}$, the operators in $\text{OP}^s$ are continuous in the sense that they map $\mathcal{H}^r$ into $\mathcal{H}^{r-s}$ for every $r\in\mathbb{R}$, and moreover, $\text{OP}^r\cdot\text{OP}^s\subseteq\text{OP}^{r+s}$. 
Additionally, if $T\in\text{OP}^r$ for some $r\geq0$, then $\delta(T)\in\text{OP}^r$. \\

Let $\nabla$ denote the derivation 
\begin{equation*}
\nabla(T):=[\Delta,T].
\end{equation*}
Under the above hypothesis, for all $T\in \mathcal{A}\cup[D,\mathcal{A}]$, for all $n\in\mathbb{N}$
\begin{equation*}
\nabla^n(T)=[\Delta,[\Delta,\ldots,[\Delta,T]]\cdots] \in\text{OP}^n
\end{equation*}
($n$ commutators), see \cite[Cor.\ B.1]{connesmoscovici}. Consider the algebra $\mathcal{DO}$ generated by the elements $\nabla^n(T)$, for all $T\in\mathcal{A}\cup[D,\mathcal{A}]$, and all $n\in\mathbb{N}$. The algebra $\mathcal{DO}$ is an analogue of the algebra of differential operators, and we have a natural filtration of $\mathcal{DO}$, $\{\mathcal{DO}^n\}_{n\in\mathbb{N}}$, given by the total power of $\nabla$ applied, i.e.\ the order of the operator denoted by $\text{ord}(\ )$.  By construction, $\mathcal{DO}$ is stable under the derivation $\nabla$ and
\begin{equation*}
\nabla(\mathcal{DO}^n)\subseteq\mathcal{DO}^{n+1}.
\end{equation*}
Also, for $A\in\mathcal{DO}^n$ and $z\in\mathbb{C}$, we have
\begin{equation*}
A\lvert D \rvert^z\in\text{OP}^{n+\text{Re}(z)}.
\end{equation*}
In \cite[Sect.\ 4.5]{higson06} it is shown that the algebra $\mathcal{DO}$ is the smallest algebra of linear operators on $\mathcal{H}^\infty$ which contains $\mathcal{A}$ and $[D,\mathcal{A}]$, and which is closed under the operation $X\mapsto\nabla(X)$, equipped with the filtration that 
\begin{enumerate}
\item If $a\in\mathcal{A}$, then $\text{ord}(a)=0$ and $\text{ord}([D,a])=0$.
\item If $X\in\mathcal{DO}$, then $\text{ord}(\nabla(X))\leq\text{ord}(X)+1$.
\end{enumerate}

Let $\mathcal{B}$ denote the algebra of operators on $\mathcal{H}^\infty$ generated by all the powers of $\delta$ of all the elements in $\mathcal{A}\cup[D,\mathcal{A}]$
\begin{equation*}
\mathcal{B}:={\text{span}}\{\delta^{n_i}(a_i),\ \delta^{n_j}([D,a_j]):\ \ n_i, n_j\in\mathbb{N}, \ a_i,a_j\in\mathcal{A}\}.
\end{equation*}
Thus, for a regular spectral triple we have that $\mathcal{B}\subseteq\text{OP}^0$, and for any $z\in\mathbb{C}$ the operator $b\lvert D \rvert^z$ lies in $\text{OP}^{\text{Re}(z)}$ for any $b\in\mathcal{B}$. Moreover, the linear span of the operators of the form $b\lvert D \rvert^k$, where $k\geq0$ is an algebra since
\begin{equation*}
b_1\lvert D \rvert^{k_1}\cdot b_2\lvert D \rvert^{k_2}=\sum_{j=0}^{k_1}\binom{k_1}{j}b_1\delta^j(b_2)\lvert D \rvert^{k_1+k_2-j},
\end{equation*}
for all $b_1,b_2\in\mathcal{B}$, and all $k_1,k_2\in\mathbb{N}$. \\

\begin{lemma}[{\cite[Lem.\ 4.27]{higson06}}] 
Assume that $(\mathcal{A},\mathcal{H},D)$ is a regular spectral triple. Every operator in $\mathcal{DO}$ of order $k$ may be written as a finite sum of operators $b\lvert D \rvert^\ell$, where $b$ belongs to the algebra $\mathcal{B}$ and $\ell\leq k$. \\
\end{lemma}

This leads to the following definition of pseudodifferential operators on a regular spectral triple: \\
\begin{definition}[{\cite[Sect.\ 2]{paycha} and \cite[Def.\ 4.28]{higson06} for the integer case}]
If $a\in\mathbb{C}$, an \emph{abstract pseudodifferential operator} of order $a$ is an operator $A:\mathcal{H}^\infty\to\mathcal{H}^\infty$, that has the following expansion:
\begin{equation}
A\sim\sum_{j=0}^\infty b_j\lvert D \rvert^{a-j},\quad b_j\in\mathcal{B},
\end{equation}
meaning by this, that for any $N\in\mathbb{N}$
\begin{equation*}
R_N(A):=A-\sum_{j=0}^{N-1}b_j\lvert D \rvert^{a-j}\in\text{\emph{OP}}^{a-N}.
\end{equation*}
\end{definition}

By \cite[Eq.\ (2.1)]{paycha} (see also \cite[Eq.\ (11)]{connesmoscovici}), for any abstract pseudodifferential operator $A$ and any complex number $\alpha$ we have
\begin{equation*}
\lvert D \rvert^\alpha A\sim A\lvert D \rvert^\alpha+\sum_{k=1}^\infty c_{\alpha,k}\delta^k(A)\lvert D \rvert^{\alpha-k},
\end{equation*}
where $c_{\alpha,k}:={\alpha(\alpha-1)\cdots(\alpha-k+1)}/{k!}$ for any positive integer $k$. This way, if $A$ has order $a$, then $\delta^k(A)\lvert D \rvert^{\alpha-k}$ has order $a+\alpha-k$, whose real value decreases as $k$ increases. \\

\begin{proposition}[{\cite[Prop.\ 4.31]{higson06}}]
The set of all pseudodifferential operators is a filtered algebra. If $T$ is a pseudodifferential operator, then so is $\delta(T)$, and moreover $\text{\emph{ord}}(\delta(T))\leq\text{\emph{ord}}(T)$. \\
\end{proposition}

In \cite[App.\ B.3]{ecksteiniochum} we can also find the study of a spectral triple with a noncommutative torus as the underlying noncommutative algebra, in order to give a description of the spectral action for it. In particular, it is proved that such a spectral triple is regular and thus, abstract pseudodifferential operators as described above can be defined on it (see also \cite[Sect.\ 3]{essouabriiochumlevysitarz}). \\

As we mentioned at the beginning of this section, in \cite{connesmoscovici} A.\ Connes and H.\ Moscovici use abstract pseudodifferential operators to prove a local index formula for arbitrary spectral triples satisfying the conditions of finite summability and simple dimension spectrum, and this is done in terms of the Dixmier trace and a residue--type extension. This was further studied in \cite{paycha} in order to define a canonical trace on certain spectral triples. \\

\

We finish this article with a few words about the calculi we have reviewed. The toroidal pseudodifferential operators given in Section \ref{sect: toroidal pdos} can be seen as a discrete version of Connes' pseudodifferential operators given in Section \ref{sect: pdos connes}. These two calculi coincide on noncommutative tori as mentioned at the end of Section \ref{sect: pdos connes}, and the operators are characterized in terms of a symbol, which is defined by using the derivations given by a torus action as expressed in Definitions \ref{def:toroidal symbol} and \ref{def:smooth toroidal symbol}. \\

The construction of the deformed pseudodifferential operators given in Section \ref{sect:deformed pdos} depends on the choice of a linear connection on the torus as a manifold, which is invariant with respect to the torus action. This allows to define an appropriate tensor calculus, which is used to define globally defined symbols (independent of the choice of coordinates) on the manifold. The resulting symbol map and quantization map are equivariant under the torus action, and this enables the author of \cite{liu} to perform some computations involving pseudodifferential operators in certain noncommutative spaces such as noncommutative tori. \\

The abstract pseudodifferential operators given in Section \ref{sect:abstractpdos} are defined in terms of the Dirac operator from the representation of a spectral triple with a noncommutative torus as the underlying noncommutative algebra, and have been used to develop geometry on that space from a spectral point of view. \\

In this article we have provided an overview of these pseudodifferential calculi on noncommutative tori, and there are still some remaining questions about them, interesting for possible future work. Here we list some of them: \\

\begin{enumerate}
\item In Section \ref{sect: pdos connes} we have shown how the toroidal pseudodifferential calculus and the Connes' pseudodifferential calculus on a fixed noncommutative torus coincide. The study of a similar relation with other pseudodifferential calculi would be useful too. \\

\item A toric manifold is a closed manifold whose diffeomorphism group contains a torus action, and a toric noncommutative manifold is obtained by the deformation of the tensor calculus over a toric manifold \cite{brainlandisujlekom}. As we have seen in Section \ref{sect:deformed pdos}, noncommutative tori are examples of toric noncommutative manifolds, on which the author of \cite{liu} has developed an appropriate pseudodifferential calculus. It would be interesting to look for pseudodifferential operators on other noncommutative manifolds, like, e.g.\ $\theta$--deformations of compact Lie groups as the ones described in \cite{wilson}. \\

\item The pseudodifferential calculi described in Sections \ref{sect: toroidal pdos}, \ref{sect: pdos connes}, \ref{sect:deformed pdos} are defined in terms of the symbols of the operators, whereas the abstract pseudodifferential operators described in Section \ref{sect:abstractpdos} depend on the Dirac operator that defines the corresponding spectral triple. It would also be interesting to study how much the algebra of abstract pseudodifferential operators on a regular spectral triple depends on the Dirac operator. \\

\item In \cite{erpyunckengroupoidpdos} there is an algebraic and geometric characterization of the classical pseudodifferential operators on a smooth manifold in terms of the tangent groupoid and its natural $\mathbb{R_+^\times}$--action. In that reference, the authors define pseudodifferential operators on arbitrary filtered manifolds, in particular the Heisenberg calculus. It would be worth investigating if the deformation along the action of a group, as described in Section \ref{sect:rieffel deformation}, can also be performed in those operators. Moreover, in \cite{cadet} there is a construction of a deformation groupoid in the framework of $C^*$--algebraic deformation quantization for toric manifolds, so one could compare it with the deformation described in Section \ref{sect:deformed pdos}.

\end{enumerate}

\bibliography{CorrectedPDOsTorusarXiv.bib}
\bibliographystyle{amsplain}

\end{document}